\documentclass[times,sort&compress,3p]{elsarticle}
\journal{arXiv}
\usepackage[labelfont=bf]{caption}
\usepackage{bm}
\usepackage{ushort}

\usepackage{amsmath,amsfonts,amssymb,amsthm,booktabs,color,epsfig,graphicx,hyperref,url}
\usepackage{dsfont}
\usepackage{diagbox}

\theoremstyle{plain}% Theorem-like structures provided by amsthm.sty
\newtheorem{theorem}{Theorem}

\newtheorem{proposition}{Proposition}
\newtheorem{lemma}{Lemma}
\newtheorem{corollary}{Corollary}

\theoremstyle{definition}

\newtheorem{remark}{Remark}
\newtheorem{example}{Example}

 \newcommand*{\bigtimes}{\mathop{\raisebox{-.5ex}{\hbox{\huge{$\times$}}}}}
\newcommand{\pRz}{\mathbb R_+}
\newcommand{\Uz}{\mathbb U}
\newcommand{\Melop}{\mathcal M}

\newcommand{\Mela}[2]{\mathcal M_{#2} [{#1}]}

\newcommand{\LpA}{\mathbb L_{\mathbb R_+}^2}
\newcommand{\Lz}{\mathbb L}
\newcommand{\Dz}{\mathbb D}
\newcommand{\Cz}{\mathbb C}
\newcommand{\Rz}{\mathbb R}

\newcommand{\Nz}{\mathbb N}
\newcommand{\Wz}{\mathbb W}
\newcommand{\1}{\mathds 1}

\newcommand{\E}{\mathbb E}
\newcommand{\Var}{\mathbb V\mathrm{ar}}

\newcommand{\nset}[1]{{\left\llbracket #1\right\rrbracket }}	
	
\newcommand{\nsetro}[1]{{\left\llbracket #1\right\llbracket }}	
\DeclareMathOperator*{\argmin}{arg\,min}
\newcommand{\rSob}{L}

\newcommand{\wSob}{s}

\DeclareMathSymbol{\minus}{\mathbin}{AMSa}{"39}

\begin{document}

\begin{frontmatter}

\title{Anisotropic spectral cut-off estimation under multiplicative measurement errors}

\author[1]{Sergio Brenner Miguel \corref{mycorrespondingauthor}}
%\author[2]{Author Two\corref{mycorrespondingauthor}}

\address[1]{Institut f\"ur Angewandte
	Mathematik, M$\Lambda$THEM$\Lambda$TIKON, Im Neuenheimer Feld 205,
	D-69120 Heidelberg, Germany}
%\address[2]{Address of Author Two in his country's language and rules}

\cortext[mycorrespondingauthor]{Corresponding author. Email address: \url{brennermiguel@math.uni-heidelberg.de}}

\begin{abstract}
We study the non-parametric estimation of an unknown density $f$ with support on $\mathbb R_+^d$ based on an i.i.d. sample with multiplicative measurement errors. 
The proposed fully-data driven procedure is based on the estimation of the Mellin
transform of the density $f$ and a regularisation of the inverse of
the Mellin transform by a spectral cut-off. The upcoming bias-variance trade-off  is dealt with by a
data-driven anisotropic choice of the cut-off parameter. In order to discuss the bias term, we consider the
\textit{Mellin-Sobolev spaces} which characterize the
regularity of the unknown density $f$ through the decay of its Mellin
transform. Additionally,  we show minimax-optimality over \textit{Mellin-Sobolev spaces} of the
spectral cut-off density estimator. 
\end{abstract}


\end{frontmatter}

\section{Introduction\label{i}}

In this work we consider  the estimation of an unknown density
$f:\mathbb R^d \rightarrow \mathbb R$ of
a positive random variable $\bm X=(X_1, \dots, X_d)$ given independent and identically
distributed (i.i.d.) copies of $\bm Y=\bm X\bm U=(X_1U_1, \dots, X_d U_d)$, where $\bm X$ and $\bm U$ are
independent of each other and $\bm U$ has a known density $g:\mathbb R^d
\rightarrow \mathbb R$. The density  $f_{\bm Y}:\mathbb R^d
\rightarrow \pRz$ of $\bm Y$ is then given by
\begin{equation*}
f_{\bm Y}(\bm y)=[f * g](\bm y):= \int_{\pRz^d} f(\bm x)g(\bm y/\bm x) \bm x^{\ushort{\minus\bm 1}}d\bm x\quad\forall \bm y\in\pRz^d,
\end{equation*}
where $\bm x/\bm y:= (x_1/y_1, \dots, x_d/y_d)$ and $\bm x^{\ushort{\minus\bm1}}:= \prod_{j=1}^d x_j^{-1}$. Here "$*$" denotes multiplicative convolution. The estimation of
$f$ using an i.i.d. sample $\bm Y_1, \dots, \bm Y_n$ from $f_{\bm Y}$ is thus an
inverse problem called
multiplicative deconvolution.  \\
In the additive deconvolution literature the density estimation for multivariate variables based on non-parametric estimators has been studied by many authors. A kernel estimator approach was investigated by \cite{ComteLacour2013} with respect to $\mathbb L^2$-risk and by \cite{Rebelles2016}  for general $\mathbb L^p$-risk. The  multivariate convolution structure density model was considered by the authors \cite{LepskiWiller2019}. The recent work \cite{Dussap2021} focuses on the study of deconvolution problems on $\pRz^d$ and introduces a data-driven estimator based on  a projection on the Laguerre basis. To the knowledge of the author,  the estimation for multivariate random variables with multiplicative measurement errors has not been studied yet.\\
For the univariate case, the recent work of \cite{Brenner-MiguelComteJohannes2020} should be mentioned which uses the Mellin transform to construct a density estimator under multiplicative measurement errors. The model of multiplicative measurement errors was motivated in the work of \cite{BelomestnyGoldenshluger2020} as a generalisation of several models, for instance the multiplicative censoring model or the stochastic volatility model. 

A summary of related work regarding the connection between the multiplicative measurement errors model and similar models can be found in \cite{Brenner-MiguelComteJohannes2020} and \cite{BelomestnyGoldenshluger2020}.
In the work of \cite{BelomestnyGoldenshluger2020}, the authors used the Mellin transform to construct a kernel estimator for the pointwise density estimation.  In their work, the authors shown that the log transformation of the observation is a special case of their estimation
strategy. In fact, by applying the logarithm the model
$Y=XU$ writes as $\log(Y)=\log(X)+\log(U)$. This naive appoach allows then the usage of commonly used deconvolution techniques to construct an estimator of the density $\log(X)$ (see for example \cite{Meister2009}) and which can be then transformed back to an estimator of $f$.
It is worth stressing out, that in this case the regularity assumptions are considered for the density of $\log(X)$ instead of $f$ direclty. This provoces difficulties for the interpretation of these regularity conditions. For the global risk case, additional complications occurs using this naive approach as pointed out by 
\cite{ComteDion2016}. \\
In this work, we generalise the results of \cite{Brenner-MiguelComteJohannes2020} 
in a similar way to the works
\cite{ComteLacour2013} and \cite{Dussap2021} for the additive deconvolution model.
To do so, we introduce a notion of the Mellin transform for multivariate random variables and show that the necessary properties of the univariate Mellin transform remain true.
Exploiting the multiplication theorem, that is $\mathcal M[f_Y]=\mathcal M[f]\mathcal M[g]$
\cite{Brenner-MiguelComteJohannes2020} introduced for the univariate case a spectral cut-off density
estimator of $f$ based on the $i.i.d.$ sample $Y_1, \dots, Y_n$.  
Considering the multivariate case, we are analogously making use of the
multiplication theorem of the Mellin transform and apply a spectral cut-off regularisation of the inversion of the Mellin-transform to define a density estimator. The accuracy of the proposed estimator
is measured in terms of the global risk with respect to a weighted
$\LpA$-norm. We identify the underlying inverse problem using the rich theory of Mellin transform and 
characterise the natural regularity
conditions expressed
in the form of Mellin-Sobolev spaces. Here, we borrow ideas from the inverse problems community
(\cite{EnglHankeNeubauer1996}) and discuss the relation between the Mellin-Sobolev spaces and analytical properties of the density $f$. In the regularisation step of the inverse problem, an additional tuning 
parameter is introduced.  For this parameter we propose a model selection method to end up with a fully data-driven estimator. We establish an
oracle inequality for the fully-data driven spectral cut-off
estimator under fairly mild assumptions on the error density
$g$. 
Moreover, we show that uniformly over Mellin-Sobolev spaces the proposed
data-driven estimator is minimax-optimal by stating both an upper and lower 
bound for the mean weighted integrated squared error of the minimax risk of the density estimation $f$ given an i.i.d. sample of $\bm Y$.

The paper is organized as follows. In Section \ref{mt} we begin with an introduction of the Mellin transform for multivariate random variables including several properties which are commonly used throughout this paper. Based on the observations $\bm Y_1, \dots, \bm Y_n$, we then introduce the spectral cut-off estimator of the density $f$ and analyse its properties for a large class of error densities. Furthermore, we study the global behavior of the proposed estimator over the Mellin-Sobolev spaces for smooth error density. Here, we show upper and lower bounds for the weighted $\LpA$-risk of our estimator implying its minimax-optimality.
In Section \ref{dd} we propose a data-driven method for the choice of the cut-off parameter only depending on the sample $\bm Y_1,\dots, \bm Y_n$ based on a model selection. 
Finally, results of a simulation
study are reported in section \ref{si} which visualize the reasonable finite
sample performance of our estimators. Proofs of theorems of Section \ref{mt} and
Section \ref{dd} are postponed to the Appendix.

\section{Minimax theory\label{mt}}

In this section we introduce the Mellin transform and collect some of its properties while more detailed proof sketches are given in Appendix \ref{a:prel}.
Define for a weight function
$\omega:\pRz^d \rightarrow \pRz$  the corresponding
weighted norm   by $\|h\|_{\omega}^2 := \int_{\pRz^d}
|h(\bm x)|^2\omega(\bm x)d\bm x $ for a measurable function $h:\pRz^d \rightarrow \Cz$. Denote by
$\Lz^2(\pRz^d,\omega)$ the set of all complex-valued, measurable functions with
finite $\|\, .\,\|_{\omega}$-norm and by $\langle h_1, h_2
\rangle_{\omega} := \int_{\pRz^d}  h_1(\bm x)
h_2(\bm x)\omega(\bm x)d\bm x$ for $h_1, h_2\in \Lz^2(\pRz^d,\omega)$
the corresponding weighted scalar product. Similarly, define $\Lz^2(\Rz^d):=\{ H:\Rz^d \rightarrow \Cz\, \text{ measurable }: \|H\|_{\Rz^d}^2:=  \int_{\Rz^d} H(\bm x)\overline{H(\bm x)} d\bm x <\infty \}$ and $\Lz^1(\Omega,\omega):=\{h: \Omega \rightarrow \Cz: \|h\|_{\Lz^1(\Omega,\omega)}:= \int_{\Omega} |h(\bm x)|\omega(\bm x)d\bm x < \infty \}$ for any $\Omega \subseteq \Rz^d$. \\
For two vectors $\bm u=(u_1,\dots, u_d)^T, \bm v=(v_1,\dots,v_d)^T \in \Rz^d$ and a scalar $\lambda \in \Rz$ we define the componentwise multiplication  $\bm u\bm v := \bm u\cdot \bm v:=(u_1v_1, \dots, u_dv_d)^T$ and denote by $\lambda \bm u$ the usual scalar multiplication. Further, if there exists no index  $i\in \nset{d}$ such that $v_i=0$  we define the multivariate power through $\bm v^{\ushort{\bm u}}:= \prod_{i\in \nset{d}} v_i^{u_i}$. Additionally, we define the componentwise division by $\frac{\bm u}{\bm v}:= \bm u /\bm v := (u_1/v_1, \dots, u_d/v_d)^T$ We denote the usual Euclidean scalar product and norm on $\Rz^d$ through $\langle \bm u, \bm v\rangle := \sum_{i\in \nset{d}} u_i v_i$ and $|\bm u|:= \sqrt{\langle \bm u, \bm u \rangle }$. 
\subsection{The Mellin transform} Let $\bm c\in \pRz^d$. For two functions $h_1,h_2\in \mathbb L^1(\pRz^d, \bm x^{\ushort{\bm c-\bm 1}})$ we define the multiplicative convolution $h_1*h_2$ of $h_1$ and $h_2$ by
\begin{align}
(h_1*h_2)(\bm y)=\int_{\pRz^d} h_1(\bm y/\bm x) h_2(\bm x) \bm x^{\ushort{\minus\bm 1}} d\bm x, \quad \bm y\in \mathbb R^d.
\end{align}
In fact, one can show that the function $h_1*h_2$ is well-defined, $h_1*h_2=h_2*h_1$ and $h_1*h_2 \in \mathbb L^1(\pRz^d,\bm x^{\ushort{\bm c-\bm1}})$. A proof sketch of this and the following properties can be found in Appendix \ref{a:prel}. Further, if additionally $h_1\in \mathbb L^2(\pRz^d, \bm x^{\ushort{2\bm c-\bm 1}})$ then $h_1*h_2 \in \mathbb L^2(\pRz^d,\bm x^{\ushort{2\bm c-\bm 1}})$. 

We will now define the Mellin transform for functions $h_1\in \mathbb L^1(\pRz^d, \bm x^{\ushort{\bm c-\bm1}})$. To do so, let $h_1\in \mathbb L^1(\pRz^d,\bm x^{\ushort{\bm c-\bm 1}})$. Then, we define the Mellin transform of $h_1$ at the development point $\bm c\in \mathbb R^d$ as the function $\mathcal M_{\bm c}[h]:\mathbb R^d\rightarrow \mathbb C$ by
\begin{align}
\mathcal M_{\bm c}[h_1](\bm t):= \int_{\pRz^d} \bm x^{\ushort{\bm c-\bm 1+i\bm t}} h_1(\bm x)d\bm x, \quad \bm t\in \mathbb R^d.
\end{align}
Note that for any density $h\in \Lz^1(\pRz^d, \bm x^{\ushort{\bm 0}})$ of a positive random variable $\bm Z$ the property $h\in \Lz^1(\pRz^d,\bm x^{\bm c-\bm 1})$ is equivalent to $\E_{h}(\bm Z^{\ushort{\bm c-\bm 1}})< \infty.$

One key property of the Mellin transform, which makes it so appealing for the use of multiplicative deconvolution, is the so-called convolution theorem, that is for $h_1, h_2\in \mathbb L^1(\pRz^d,\bm x^{\ushort{\bm c-\bm 1}})$ holds
\begin{align}
\mathcal M_{\bm c}[h_1*h_2](\bm t)=\mathcal M_c[h_1](\bm t) \mathcal M_c[h_2](\bm t), \quad \bm t\in \mathbb R^d.
\end{align}
In analogy to the Fourier transform, one can define the Mellin transform for functions $h\in \Lz^2(\pRz^d, \bm x^{\ushort{2\bm c-\bm1}})$. In fact, 
let $\varphi:\Rz^d \rightarrow \pRz^d, \bm x\mapsto (\exp(-x_1), \dots, \exp(-x_d))^T$ and $\varphi^{-1}: \pRz^d \rightarrow \Rz^d$ its inverse. Then as diffeomorphisms $\varphi$ and $ \varphi^{-1}$ map  Lebesgue null sets on Lebesgue null sets. Thus the isomorphism $\Phi_{\bm c}:\Lz^2(\pRz^d,\bm x^{\ushort{2\bm c-\bm 1}}) \rightarrow \Lz^2(\Rz^d), h\mapsto \varphi^{\bm c} h\circ \varphi$ is well-defined for any $\bm c \in \Rz^d$. Furthermore, let $\Phi^{-1}_{\bm c}: \Lz^2(\Rz^d) \rightarrow \Lz^2(\pRz^d, \bm x^{\ushort{2\bm c-\bm 1}})$ denote its inverse. Then for $h\in \Lz^2(\pRz^d,\bm x^{\ushort{2\bm c-\bm 1}})$ we define the Mellin transform of $h$ developed in $\bm c\in \Rz^d$ by
\begin{align}\label{eq:mel:def}
\Mela{h}{\bm c}(\bm t):= (2\pi)^{d/2} \mathcal F[\Phi_{\bm c}[h]](\bm t), \quad \bm t\in \Rz^d,
\end{align} 
where $\mathcal F: \Lz^2(\Rz^d)\rightarrow \Lz^2(\Rz^d)$ is the Plancherel-Fourier transform. Due to this definition several properties of the Mellin transform can be deduced from the well-known theory of Fourier transforms. In the case $h\in \Lz^1(\pRz^d,\bm x^{\ushort{\bm c-\bm 1}}) \cap \Lz^2(\pRz^d,\bm x^{\ushort{2\bm c-\bm 1}})$ we have 

\begin{align}\label{eq:mel:l1}
\Mela{h}{\bm c}(\bm t) =\int_{\pRz^d} \bm x^{\ushort{\bm c-\bm 1+i\bm t}} h(\bm x)d\bm x, \quad \bm t\in \Rz^d,
\end{align}
which coincides with the usual notion of Mellin transforms as considered in \cite{ParisKaminski2001} for the case $d=1$. \\
Further, due to this construction of the operator $\mathcal M_{\bm c}: \Lz^2(\pRz^d,\bm x^{\ushort{2\bm c-\bm 1}}) \rightarrow \Lz^2(\Rz^d)$ it is an isomorphism and we denote by $\mathcal M_{\bm c}^{-1}: \Lz^2(\Rz^d) \rightarrow \mathcal \Lz^2(\pRz^d,\bm x^{\ushort{2\bm c-\bm 1}})$ its inverse. If additionally to $H\in \Lz^2(\Rz^d)$, $H\in \Lz^1(\Rz^d)$ holds then we can express the inverse Mellin transform explicitly through
\begin{align}\label{eq:Mel:inv}
\Melop_{\bm c}^{-1}[H](\bm x)= \frac{1}{(2\pi)^d } \int_{\Rz^d} \bm x^{\ushort{-\bm c-i\bm t}} H(\bm t) d\bm t, \quad \text{ for any } \bm x\in \pRz^d.
\end{align} 
Furthermore, we can directly show that a Plancherel-type equation holds for the Mellin transform, that is for all $ h_1, h_2 \in \Lz^2(\pRz^d,\bm x^{\ushort{2\bm c-\bm 1}})$ holds 
\begin{align}\label{eq:Mel:plan}
\hspace*{-1cm}\langle h_1, h_2 \rangle_{\bm x^{\ushort{2\bm c-\bm 1}}} = \frac{1}{(2\pi)^{d}} \langle \Mela{h_1}{\bm c}, \Mela{h_2}{\bm c} \rangle_{\Rz^d}\, \, \text{ and thus } \, \, \| h_1\|_{\bm x^{\ushort{2\bm c-\bm 1}}}^2=\frac{1}{(2\pi)^{d}}  \|\Mela{h_1}{\bm c}\|_{\Rz^d}^2.
\end{align} 

\subsection{Estimation strategy}

Let us define for $k\in \pRz^d$ the hyper cuboid $Q_{\bm k}:=\{\bm x \in \Rz^d: \forall i\in \nset{d}: |\bm x_i |\leq k_i\}$. Then for $f\in \Lz^2(\pRz^d, \bm x^{\ushort{2\bm c-\bm1}})$ we have that $\Mela{f}{\bm  c} \mathds 1_{Q_{\bm k}} \in \Lz^2(\Rz^d)\cap \Lz^1(\Rz^d)$ and thus
\begin{align*}
f_{\bm k}(\bm x):= \frac{1}{(2\pi)^d} \int_{Q_{\bm k}} \bm x^{\ushort{-\bm c-i\bm t }}\Mela{f}{\bm c}(\bm t) d\bm t, \quad \bm x \in \pRz^d, 
\end{align*}
is an approximation of $f$ in the $\Lz^2(\pRz^d, \bm x^{\ushort{2\bm c-\bm1}})$-sense, that is $\|f_{\bm k}-f\|_{\bm x^{2\bm c-\bm1}} \rightarrow 0$ for $\bm k\rightarrow \bm \infty$ where the limit $\bm k\rightarrow \bm  \infty$ means that every component of $\bm k$ is going to infinity. 

Now let us additionally assume that $f\in \Lz^2(\pRz^d, \bm x^{\ushort{2\bm c-\bm1}}) \cap \Lz^1(\pRz^d, \bm x^{\ushort{\bm c-\bm1}})$ and $g\in \Lz^1(\pRz^d, \bm x^{\ushort{\bm c-\bm 1}})$. Then from the convolution theorem, one deduces  $\Mela{f_{\bm Y}}{\bm c}= \Mela{f}{\bm c}\Mela{g}{\bm c}$. Under the mild assumption that $\Mela{g}{\bm c}(t) \neq 0$, for any $t\in \Rz^d$ we can rewrite the last equation as $\Mela{f}{\bm c}= \Mela{f_{\bm Y}}{\bm c}/\Mela{g}{\bm c}$. Thus we have 
\begin{align*}
f_{\bm k}(\bm x):= \frac{1}{(2\pi)^d} \int_{Q_{\bm k}} \bm x^{\ushort{-\bm c-i\bm t}}\frac{\Mela{f_{\bm Y}}{\bm c}(\bm t)}{\Mela{g}{\bm c}(\bm t)} d\bm t, \quad\bm x \in \pRz^d. 
\end{align*} 
Let us now consider for any $\bm t\in \Rz^d$ the unbiased estimator $\widehat{\mathcal M}_c(\bm t):= n^{-1} \sum_{j\in \nset{n}} \bm Y_j^{\ushort{\bm c-\bm 1+i\bm t}}$ of $\Mela{f_{\bm Y}}{\bm c}(\bm t)$. We see easily that $|\widehat{\mathcal M}_c(\bm t)| \leq |\widehat{\mathcal M}_c(\bm 0)|< \infty$ almost surely. If additionally $\1_{Q_{\bm k}} \Mela{g}{\bm c}^{-1} \in \Lz^2(\Rz^d)$ then $\1_{Q_{\bm k}} \widehat{\mathcal M}_{\bm c}/\Mela{g}{\bm c} \in \Lz^2(\Rz^d) \cap \Lz^1(\Rz^d)$ and we can define our spectral cut-off density estimator by $\widehat f_{\bm k} :=\mathcal M^{-1}_{\bm c} [\1_{Q_{\bm k}} \widehat{\mathcal M}_{\bm c}/ \Mela{g}{\bm c}]$. More explicitly, we have 
\begin{align}\label{eq:est:def}
\widehat f_{\bm k}(\bm x)= \frac{1}{(2\pi)^d} \int_{Q_{\bm k}} \bm x^{\ushort{-\bm c-i\bm t}} \frac{\widehat{\mathcal M}_{\bm c}(\bm t)}{\Mela{g}{\bm c}(\bm t)} d\bm t,\quad  \bm x \in \pRz^d.
\end{align}
Up to now, we had two minor assumptions on the error density $g$ which we want to collect in the following assumption:
\begin{align} \label{eq:ass:g0}
\forall \bm t \in \Rz^d: \Mela{g}{\bm c}(\bm t)\neq 0 \text{ and } \forall \bm k \in \pRz^d: \int_{Q_{\bm k}} |\Mela{g}{\bm c}(\bm t)|^{-2} d\bm t <\infty.
\end{align}
\begin{remark}
	The assumption \eqref{eq:ass:g0} resembles strongly the  rather typical error assumption in context of deconvolution problems, compare \cite{Meister2009}. Examples of multivariate density, which fullfils the assumption \eqref{eq:ass:g0} are given in Example \ref{ex:lognorm} and Example\ref{ex:unif}. It is worth stressing out, that one can construct deconvolution estimators under a weaker assumption on the error density, that is $\mathcal M_{\bm c}[g](\bm t) \neq 0$ almost every, compare \cite{Brenner-Miguel2021}. In particular, the weaker assumption is in that sense minimal that without it, we need addtionall constraints on the class of densities to ensure that the density is indeed identifiable, compare \cite{Meister2009} and \cite{Brenner-Miguel2021} .
	\end{remark}
The following proposition shows that the proposed estimator is consistent for a suitable choice of the cut-off parameter $\bm k \in \pRz^d$. Its proof is postponed to Appendix \ref{a:mt}.

\begin{proposition}\label{pr:consis}
	Let $f\in \Lz^2(\pRz^d, \bm x^{\ushort{2\bm c-\bm 1}})$, $\sigma:= \E_{f_Y}(\bm Y^{\ushort{2\bm c- \bm 2}})< \infty$ and assume that \eqref{eq:ass:g0} holds for $g$. Then we have for any $\bm k\in \pRz^d$,
	\begin{align} \label{eq:consis}
	\E_{f_Y}^n(\|f-\widehat f_{\bm k}\|_{\bm x^{\ushort{2\bm c-\bm 1}}}^2) = \|f-f_{\bm k}\|_{\bm x^{\ushort{2\bm c-\bm 1}}}^2+ \frac{\sigma \Delta_g(\bm k)}{n},
	\end{align}
	where $\Delta_g(\bm k):= \frac{1}{(2\pi)^d}\int_{Q_{\bm k}} |\Mela{g}{\bm c}(\bm t)|^{-2} d\bm t$. Now choosing $\bm k_n$ such that $\Delta_g(\bm{k_n})n^{-1} \rightarrow 0$ and $\bm k_n \rightarrow \bm \infty$ implies the consistency of $\widehat f_{\bm{k_n}}$.
\end{proposition}
Let us comment on the last result. For a suitable choice of the spectral cut-off parameter $\bm k\in \pRz^d$ we can show that the estimator is consistent in the sense of the weighted $\Lz^2$ distance. The second parameter, the model parameter $\bm c\in \Rz^d$, is linked to the considered risk and the assumptions on the densities $f$ and $g$.  In fact, choosing $\bm c = \bm 1$, we see that $\sigma = 1$ for any densities $f,g$. In this scenario, no additional moment assumptions on the densities $f$ and $g$ are needed. If one wants to consider the unweighted $\mathbb L^2$, one should set $\bm c=\bm{1/2}$ which leads in the case $d=1$ to the assumption that $\mathbb E_{f_Y}(Y^{-1})< \infty$. As one can see, the choice of the parameter $c\in \Rz$ is more of a modeling nature. Nevertheless, it is worth stressing out, that the authors of \cite{BelomestnyGoldenshluger2020} considered optimal choices of $c\in \Rz$ for the pointwise estimation of the density. Since in the global estimation the risk itself is dependent on $c\in \Rz$, the role of $c\in \Rz$ is quite different from the pointwise estimation.

Up to now, the assumptions on $f$ and $g$ were to ensure the well-definedness of the estimator and the weighted $\Lz^2$-risk. Here, we can already see that the first summand, called bias term, in Proposition \ref{pr:consis} is decreasing if $\bm k\in \pRz$ is increasing in any direction while the second summand, called variance term, is increasing. For a more sophisticated analysis of both terms we will consider stronger assumptions on the densities $f$ and $g$. Let us first start with the noise density $g$.
\subsection{Noise assumption}
As already mentioned, the variance term in \eqref{eq:consis} is monotonically increasing in each component $k_j, j\in \nset{d},$ of $\bm k$.  More precisely, the growth of $\Delta_g$ is determined by the decay of the Mellin transform of $g$ in each direction. \\
In the context of additive deconvolution problems (compare \cite{Fan1991}), densities whose Fourier transform decay polynomially, like in Examples \ref{ex:beta} and \ref{ex:log-gamma}, are called smooth error densities. To stay in this way of speaking we say that an error density $g$ is a smooth error density if there exists $c_g,C_g\in \pRz$ such that
\begin{align}\label{eq:ass:g1}
 c_g \prod_{j\in \nset{d}} (1+t_j^2)^{-\gamma_j/2} \leq |\Mela{g}{\bm c}(\bm t)| \leq C_g \prod_{j\in \nset{d}} (1+t_j^2)^{-\gamma_j/2}, t\in \mathbb R.
\end{align} 
This assumption on the error density  was also considered in the works of \cite{BelomestnyGoldenshluger2020} and \cite{Brenner-MiguelComteJohannes2020}. 
Under this assumption, we see that $\Delta_g(\bm k)\leq C_g \prod_{j\in \nset{d}} k_j^{2\gamma_j +1}$ for every  $\bm k \in \pRz^d$. 
After a more sophisticated bound of the variance term we will consider now the bias term which occurs in  \eqref{eq:consis}.
\subsection{Regularity spaces}
Let us for $\bm s, \bm c \in \pRz^d$ define the anisotropic Mellin-Sobolev space by
\begin{align}\label{eq:ani:mell:sob}
\Wz^{\bm s}_{\bm c}(\pRz^d):= \{ h\in \Lz^2(\pRz, \bm x^{\ushort{2\bm c-\bm 1}}): |h|_{\bm s, \bm c}^2:= \sum_{j\in \nset{d}}\| (1+t_j)^{s_j} \Mela{h}{\bm c}\|_{\Rz^d}^2< \infty \}
\end{align}
and the corresponding ellipsoids with $L\in \pRz$ by $\Wz^{\bm s}_{\bm c}(L):=\{ h\in \Wz^{\bm s}_{\bm c}(\pRz^d): |h|_{\bm s, \bm c}^2 \leq L\}$.
Since $\bigcup_{j=1}^d\{ \bm t \in \Rz^d: |t_j|>k_j\} \supset Q_{\bm k}^c$ we deduce from the assumption $f\in \Wz^{\bm s}_{\bm c}(L)$ that
\begin{align*}
\|f-f_{\bm k}\|_{\bm x^{2\bm c- \bm 1}}^2 \leq \sum_{j\in \nset{d}} \| \1_{[-k_j, k_j]^c} \Mela{f}{\bm c} \|_{\Rz^d}^2 \leq L \sum_{j\in \nset{d}} k_j^{-2s_j}.
\end{align*}
Setting $\Dz^{\bm s}_{\bm c}(L):=\{f \in \Wz^{\bm s}_{\bm c}(L): f \text{ density}, \E_{f}(\bm X^{\ushort{2\bm c-\bm 2}}) \leq L)$, the previous discussion leads to the following statement.
\begin{lemma}\label{lem:upper:minimax}
	Let $f\in \Lz^2(\pRz^d, \bm x^{\underline{2\bm c- \bm1}})$ and $\E_{g}(\bm U^{\ushort{2\bm c-\bm 2}})< \infty.$ Then under the assumptions \eqref{eq:ass:g0} and \eqref{eq:ass:g1} holds
	\begin{align*}
	\sup_{f\in \Dz^{\bm s}_{\bm c}(L)} \E_{f_Y}^n(\|f- \widehat f_{\bm k_n}\|_{\bm x^{\ushort{2\bm c-\bm 1}}})^2 \leq C(L,g ,\bm s) n^{-1/(1+2^{-1}\sum_{j\in \nset{d}} (2\gamma_j+1) s_j^{-1})}
	\end{align*}
	for the choice $\bm k_n=(k_{1,n}, \dots, k_{d,n})$ with $k_{i,n}:= n^{1/(2s_i+s_i\sum_{j\in \nset{d}} (2\gamma_j+1)s_j^{-1})}$.
\end{lemma}

Considering the rate in Lemma \ref{lem:upper:minimax} the natural question  arises if whether exists an estimator based on the sample $Y_1, \dots, Y_n$ which has a sharper rate uniformly over $\mathbb D_{\bm c}^{\bm s}(L)$. 

In the following paragraph we will show that such a scenario cannot occur. From this we deduce that our estimator $\widehat f_{k_n}$ is minimax-optimal over the ellipsoids $\mathbb D_{\bm c}^{\bm s}(L)$ for many classes of error densities.
\subsection{Lower bound}
For the following part, we will need to have further assumption on the error density $g$. In fact, we will distinguish if $c_j\in (0,1/2]$ or $c_j > 1/2$ for $j\in \nset{d}$. Let therefore $\tilde{\bm c}:= (\tilde c_1, \dots, \tilde c_d) \in \pRz^d$ where $\tilde c_j = 2^{-1} \1_{(1/2, \infty)}(c_j)$. Let us assume that $g$ has a bounded support, that is for all
$\bm x \in \pRz^d\setminus (0,1)^d: g(\bm x)=0$ and that there exists constants $c_g, C_g>0$ such that
\begin{align}\label{eq:ass:g2}
 c_g \prod_{j\in \nset{d}} (1+t_j^2)^{-\gamma_j/2} \leq |\Mela{g}{\bm{\tilde c}}(\bm t)| \leq C_g \prod_{j\in \nset{d}} (1+t_j^2)^{-\gamma_j/2},\quad  \bm t\in \mathbb R^d.
\end{align} 
With this additional assumption we can show the following theorem where its proof can be found in  Appendix \ref{a:mt}.

\begin{theorem}\label{theorem:lower_bound}
	Let $\bm \wSob, \bm \gamma \in\Nz^d$ , $\bm c>\bm{0}$ and assume that \eqref{eq:ass:g1} and \eqref{eq:ass:g2} holds. Then there exist
	constants $C_{g},L_{\bm\wSob,g,\bm c}>0,$ such that for all
	$L\geq L_{\bm \wSob,g,\bm c}$, $n\in \Nz$ and for any estimator $\widehat f$ of $f$ based
	on an i.i.d. sample $(Y_i)_{i\in \nset{n}}$, 
	\begin{align*}
	\sup_{f\in\mathbb D^{\bm s}_{\bm c}(L)}\E_{f_Y}^n(\|\widehat f-f\|_{\bm x^{\ushort{2\bm c-\bm 1}}}^2) \geq c_{g,\bm c}   n^{-1/(1+2^{-1}\sum_{j\in \nset{d}} (2\gamma_j+1) s_j^{-1})}.
	\end{align*}
\end{theorem}

Let us shortly comment on the assumption \eqref{eq:ass:g2}. For the case that $c_j>1/2$ for $j\in \nset{d}$ then $\tilde c_j = 1/2$, thus we need to assume that $\E_g(U_j^{-1/2})< \infty$ which is a mild condition. If $c_j\leq 1/2$ we have that $\tilde c_j = 0$. In this case, we have automatically that $\E_g(U_j^{-1}) <\infty$ if we assume that $\E_g(U_j^{2c_j-2})< \infty$, compare Proposition \ref{pr:consis}.

\section{Data-driven method}\label{dd}
Allthough we have shown that in certain situations the estimator $\widehat f_{\bm k_n}$ in Lemma \ref{lem:upper:minimax} is minimax-optimal, the choice of $\bm k_n$ is still dependent on the regularity parameter $\bm  s\in \pRz^d$ of the unknown density $f$, which is again, unknown. Therefore, we will propose a fully data-driven choice of $\bm k\in \pRz^d$ based on the sample $Y_1, \dots, Y_n$. 
For the special case of $d=1$ the authors of \cite{Brenner-MiguelComteJohannes2020} proposed a data-driven choice for the parameter $k\in \pRz$ based on a penalized contrast approach. For the multivariate case, a model selection approach has been mainly used if one considers an isotropic choice of the cut-off parameter, that is, instead of considering the estimator defined in \eqref{eq:est:def} one would use for $k\in \pRz$ and $B_k:=\{\bm x\in \Rz^d: |\bm x| \leq k\}$ the estimator
\begin{align*}
\widetilde f_{k}(\bm x)= \frac{1}{(2\pi)^d} \int_{B_{ k}} \bm x^{\ushort{-\bm c-i\bm t}} \frac{\widehat{\mathcal M}_{\bm c}(\bm t)}{\Mela{g}{\bm c}(\bm t)} d\bm t,\quad  \bm x \in \pRz^d.
\end{align*}
For the family $(\widetilde f_{k})_{k\in \pRz}$ a data-driven choice of the parameter $k\in \pRz$ based on a model selection approach is possible. Although it might be tempting to use this estimator as the multivariate generalisation of the estimator presented in \cite{Brenner-MiguelComteJohannes2020}, an anisotropic estimator has the advantage that it is more flexibel. In fact, if the regularity in two directions  of the density differs substantial, respectively the decay of the Mellin transform of the error density, an isotropic choice of the cut-off parameter is obviously inappropiate. For the anisotropic estimator defined in \eqref{eq:est:def} we propose a data-driven choice based on a model selection which can be used even for anistropic choices of the cut-off parameter $\bm k\in \pRz^d$. To the knowledge of the authors the usage of a model selection approach instead of a Lepski approach, compare \cite{Dussap2021} and \cite{ComteLacour2013}, has not been considered so far.
Let us reduce the set of possible parameters to 
\begin{align*}
\mathcal K_n:=\{\bm k\in \mathbb N^d: \Delta_g(\bm k) \leq n\}
\end{align*}
and define for $\bm k\in \pRz^d$ and $\chi>0$ the penalty term \begin{align*} 
\mathrm{pen}(\bm k):= \chi  \sigma \Delta_{g}(\bm k)n^{-1}.
\end{align*}
It can be seen that the bias $\|f-f_{\bm k}\|_{\bm x^{\ushort{2\bm c-\bm 1}}}^2= \|f\|_{\bm x^{\ushort{2\bm c-\bm 1}}}^2 -  \|f_{\bm k}\|_{\bm x^{\ushort{2\bm c-\bm 1}}}^2$ behaves like $- \|f_{\bm k}\|_{\bm x^{\ushort{2\bm c-\bm 1}}}^2$. Exchanging  $-\|f_{\bm k}\|_{\bm x^{\ushort{2\bm c-\bm 1}}}^2$ and $\mathrm{pen}(\bm k)$ with their empirical counterparts, we define the fully data-driven model selection by
\begin{align}\label{eq:data:driven}
\widehat{ \bm k}:= \argmin_{\bm k\in \mathcal K_n} (-\|\widehat f_{\bm k}\|_{\bm x^{\ushort{2\bm c-\bm 1}}}^2+ \widehat{\mathrm{pen}}(\bm k)), \quad \widehat{\mathrm{pen}}(\bm k):= \chi \widehat\sigma \Delta_g(\bm k)n^{-1}
\end{align}
where $\widehat{\sigma}:=n^{-1}\sum_{j\in \nset{n}} \bm Y_j^{\ushort{2\bm c-\bm 2}}$.
Now, let us show that this data-driven procedure mimics the optimal choice up to a neglectible term.
\begin{theorem}\label{theo:adap:aniso}
	Assume that $\E_{f_{\bm Y}}(\bm Y^{\ushort{7(\bm c-\bm 1)}})<\infty$, $\|\bm x^{\ushort{2\bm c-\bm1}} f_{\bm Y}\|_{\infty}<\infty$ and \eqref{eq:ass:g1} is fulfilled.
	Then for $\chi \geq 144$ we have
	\begin{align*}
	\E_{f_{\bm Y}}^n(\|f-\widehat f_{\widehat k}\|_{\bm x^{\ushort{2\bm c-\bm 1}}}^2) \leq 3 \inf_{\bm k\in \mathcal K_n} (\|f-f_{\bm k}\|_{\bm x^{\ushort{2\bm c-\bm 1}}}^2+\mathrm{pen}(\bm k) )+ \frac{C_2}{n}
	\end{align*}
	 $C_2>0$ is a constant depending on $\chi,\|f_{\bm Y}\bm x^{\ushort{2\bm c-\bm1}
	}\|_{\infty},\sigma$,$ \E_{f_{\bm Y}}(\bm Y^{\ushort{7(\bm c-\bm 1)}})$ and $ g$.
\end{theorem}
For every $\bm s\in \pRz^d$ we can see that $\bm k_n$, defined in Lemma \ref{lem:upper:minimax} lies in $\mathcal K_n$. Due to this, and the consideration in the minimax theory section, we can deduce the following Corollary directly whose proof is thus omitted.

\begin{corollary}\label{lem:dd:mm}
	Under the assumption of Theorem \ref{theo:adap:aniso} and the additional assumption that $f\in \mathbb D_{\bm c}^{\bm s}(L)$ we get
	\begin{align*}
	\E_{f_{\bm Y}}^n(\|f-\widehat f_{\widehat k}\|_{\bm x^{2\bm c-\bm 1}}^2) \leq C n^{-1/(1+2^{-1} \sum_{j\in \nset{d}} (2\gamma_j+1) s_j^{-1})}
	\end{align*}
	where $C$ is a positive constant depending on $\chi, \|f_{\bm Y} \bm x^{\ushort{2\bm c-\bm 1}}\|_{\infty}$, $\sigma$, $\E_{f_{\bm Y}}(\bm Y^{\ushort{7(\bm c-\bm1)}})$, $g$ and $L$.
\end{corollary}

\section{Examples and Numerical study}\label{si}
\subsection{Examples}
In this subsection, we aim to motivate the definition of the Mellin-Sobolev spaces and the noise assumption by considering various examples presented in work \cite{Brenner-MiguelComteJohannes2020}. For the sake of readability, we will begin with the case $d=1$ and will then consider then consider examples for $d=2$.
\paragraph{Univariate case}
\begin{example}[The Beta and the Log-Gamma Distribution] \label{ex:beta}
	Consider the family $(g_b)_{b\in \Nz}$ of $\mathrm{Beta}_{(1, b)}$-densities given by
	$
	g_b(x):= \1_{(0,1)}(x) b(1-x)^{b-1},\text{ for } b\in \mathbb N \text{ and } x\in \pRz$ 
	Obviously, we see that $g_b\in \Lz^2_{\Rz^+}(x^{2c-1})\cap \Lz^1_{\Rz^+}(x^{c-1})$  for $c>0$ and it holds
	\begin{align*}
	\Mela{g_b}{c}(t) =\prod_{j=1}^{b} \frac{j}{c-1+j+it}, \quad t\in \Rz.
	\end{align*} 
	Considering the decay of the Mellin transform we get
	$
	c_{g,c} (1+t^2)^{-b/2} \leq	|\mathcal M_c[g_b](t)| \leq C_{g,c} (1+t^2)^{-b/2},t\in \Rz,
	$
	where $c_{g,c}, C_{g, c} >0$ are positive constants only depending on $g$ and $c$.
\end{example}
\begin{example}[The Scaled Log-Gamma Distribution]\label{ex:log-gamma}
	Consider the family $(g_{\mu, a, \lambda})_{(\mu,a,\lambda) \in \Rz\times \Rz^+\times \Rz^+}$ of $\mathrm{sL}\Gamma_{(\mu, a, \lambda)}$  densities with 
	$
	g_{\mu, a, \lambda}(x)=\frac{\exp(\lambda \mu)}{\Gamma(a)} x^{-\lambda-1} (\log(x)-\mu)^{a-1}\1_{(e^{\mu}, \infty)}(x), \text{ for } a,\lambda, x \in \Rz^+ \text{ and } \mu \in \Rz.
	$Then for $c<\lambda+1$ holds $g_{\mu,a,\lambda }\in \Lz^2_{\Rz^+}(x^{2c-1})\cap \Lz^1_{\Rz^+}(x^{c-1})$ and
	\begin{align*}
	\Mela{g_{\mu, a, \lambda}}{c}(t)= \exp(\mu(c-1+it)) (\lambda-c+1-it)^{-a}, \quad t\in \Rz.
	\end{align*}
	If $a=1$ then $g_{\mu,1,\lambda}$ is the density of a Pareto distribution with parameter $e^{\mu}$ and $\lambda$. If $\mu=0$ we have that $g_{0, a, \lambda}$ is the density of a Log-Gamma distribution. 
	Considering the decay of the Mellin transform we get
	$
	c_{g,c} (1+t^2)^{-a/2} \leq	|\mathcal M_c[g_{\mu,a,\lambda}](t)| \leq C_{g,c} (1+t^2)^{-a/2}, t\in \Rz
	$,
	where $c_{g,c}, C_{g, c} >0$ are positive constants only depending on $g$ and $c$.
\end{example}
Example \ref{ex:beta} and Example \ref{ex:log-gamma} both fulfill the noise assumption \eqref{eq:ass:g1}. Furthermore, we see that their Mellin transforms are everywhere non-zero. 
In the context of the Mellin-Sobolev spaces defined in \eqref{eq:ani:mell:sob}, we see that for the density of a Beta Distribution $g_{b}\in \mathbb W_c^s(\pRz)$ if $s< b-1/2$ and analogously we can see for the density of a scaled Log-Gamma distribution that $g_{\mu, a, \lambda} \in \mathbb W^s_c(\Rz),$ if $s< a-1/2$.
Let us now consider three examples of densities with an exponential decay of the corresponding Mellin transform.
\begin{example}[Gamma Distribution]\label{ex:gam}
	Consider the family $(g_d)_{d\in \Rz_+}$ of $\Gamma_{(d, 1)}$ densities with
	$
	g_d(x) = \frac{x^{d-1}}{\Gamma(d)} \exp(-x) \1_{\Rz^+}(x)$ for $d,x\in \Rz^+.
	$
	Obviously, we see that  $g_d\in \Lz^2_{\Rz^+}(x^{2c-1})\cap \Lz^1_{\Rz^+}(x^{c-1})$ for $c>-d+1$ and it holds
	\begin{align*}
	\mathcal M_c[g_d](t)= \frac{\Gamma(c+d-1+it)}{\Gamma(d)}, \quad t\in \Rz.
	\end{align*}
	Applying the Stirling formula, compare with \cite{BelomestnyGoldenshluger2020}, leads to
	$
	c_{c,d} |t|^{2c+2d-3}\exp(-\pi |t|)\leq |\mathcal M_c[g_d](t) |^2\leq C_{c, d} |t|^{2c+2d-3}\exp(-\pi |t|),$ $ |t|\rightarrow \infty.
	$
\end{example}
\begin{example}[Weibull Distribution]\label{ex:wei}
	Consider the family $(g_m)_{m\in \pRz}$ of $\mathrm W_{(m, 1)}$ densities
	$
	g_m(x) = m x^{m-1} \exp(-x^m) \1_{\pRz}(x)$ for $m,x\in \pRz.
	$Obviously, we see that $\mathcal M_c[g_m]$ is well-defined for $c>-m+1$ and it holds
	\begin{align*}
	\mathcal M_c[g_m](t)= \frac{(c-1+it)}{m}\Gamma\left(\frac{c-1+it}{m}\right), \quad t\in \Rz.
	\end{align*}
	Applying the Stirling formula one sees that
	$
	c_{c, m} |t|^{\frac{2(c-1)+m}{m}}e^{-\pi |t|/m} \leq |\mathcal M_c[g_m](t) |^2\leq C_{c, m} |t|^{\frac{2(c-1)+m}{m}}e^{-\pi |t|/m},$   for all $|t|\rightarrow \infty.
	$
\end{example}

Obviously the Example \ref{ex:gam}, Example \ref{ex:wei}, do not satisfy assumption \eqref{eq:ass:g1}. In fact, they are examples for so-called \text{super smooth densities} which are characterized by an exponential decay of their Mellin transform, compare \cite{BelomestnyGoldenshluger2020} and \cite{Brenner-MiguelComteJohannes2020}. Although, obviously $g_d, g_m, g_{\mu, \lambda}\in \mathbb W^s_c(\pRz)$ for all choices $s\in \pRz$, the decay of the bias used in Corollary \ref{lem:upper:minimax} is more of a pessimistic nature. For these three examples, we can show that the bias is of exponential decay rather than of polynomial decay. The following table shows the resulting risk rates combining the  Example \ref{ex:beta} to Example \ref{ex:wei}. \\
\begin{table}[ht]
	\caption{Comparison of different combinations of unknown densities $f$ and error densities $g$. The decay of the Mellin transform of $f$ is in both cases polynomial. For Examples \ref{ex:beta} and \ref{ex:log-gamma} the decay of the Mellin transform of $g$ is polynomial and exponential for the Examples \ref{ex:gam} and \ref{ex:wei}.}
\hrule
	\centerline{\begin{tabular}{@{}cc|cccc@{}}
			&& $g$ & & & \\
			&\hspace*{0.6cm}$f$& $\mathrm{Beta}_{(1, b_2)}$ & $\mathrm{sL}\Gamma_{(0, a_2, 1)}$& $\Gamma_{(d_2, 1)}$ & $\mathrm W_{(m_2, 1)}$\\\midrule
			& $\mathrm{Beta}_{(1, b_1)}$& $n^{-\frac{2b_1-1}{2(b_1-b_2)}}$&$n^{-\frac{2b_1-1}{2(b_1-a_2)}}$ & $\log(n)^{1-2b_1}$&$ \log(n)^{1-2b_1} $\\
			&$\mathrm{sL}\Gamma_{(0, a_1, 1)}$& $n^{-\frac{2a_1-1}{2(a_1-b_2)}}$& $n^{-\frac{2a_1-1}{2(a_1-a_2)}}$& $\log(n)^{1-2a_1}$& $\log(n)^{1-2a_1}$  \\
	\end{tabular}}
\hrule
	\label{table:num_res}
\end{table}
We want to stress out, that there are several assumptions on the upcoming parameters of the density, for instance, Log-Gamma Distributions are only weighted square-integrable if $a>1/2$, or their interplay with the model parameter $c\in \Rz$, that we will not mention for the sake of readability of the presented table.
Considering the Examples \ref{ex:gam} and \ref{ex:wei}  for the unknown density and with smooth error density given in the Examples \ref{ex:beta} and \ref{ex:log-gamma}.
\begin{table}[ht]
		\caption{Comparison of different combinations of unknown densities $f$ and error densities $g$. The decay of the Mellin transform of $f$ is in all three cases exponential. The decay of the Mellin transform of $g$ is in both cases polynomial.  \label{table:num_res2}}
	\hrule
	\centerline{\begin{tabular}{@{}cc|cc@{}}
			&&$g$& \\
			&\hspace*{0.35cm}$f$& $\mathrm{Beta}_{(1, b_2)}$ & $\mathrm{sL}\Gamma_{(0, a_2, 1)}$\\\midrule
			& $\Gamma_{(d_1, 1)}$  & $\log(n)^{2b_2+1} n^{-1}$&$\log(n)^{2a_2+1} n^{-1}$ \\
			& $\mathrm W_{(m_1, 1)}$& $\log(n)^{2b_2+1} n^{-1}$&$\log(n)^{2a_2+1} n^{-1}$  \\
	\end{tabular}}

	\hrule
\end{table}\vspace*{1cm}

We will now consider two examples of densities $g$ which do not factorise, that is, where there exists no $g_1, g_2$ with $g(\bm x)= g_1(x_1)g_2(x_2)$ since in this case, the Mellin transform $\mathcal M_{\bm c}[g](\bm t)= \mathcal M_{c_1}[g_1](t_2) \mathcal M_c[g_2](t_2)$.

\begin{example}[Bivariate Log-Normal Distribution]\label{ex:lognorm}
	Consider the family $(g_{\bm\mu,\bm\Sigma})_{(\bm\mu, \bm \Sigma) \in \Rz^2\times \pRz^{(2,2)}}$ of $\mathrm{LN}_{(\bm \mu, \bm\Sigma)}$ densities where $g_{\bm \mu,\bm \Sigma}$ for $\bm \Sigma\in \pRz^{(2,2)}$, positive definit, $\bm x \in \pRz^2$ and $\bm\mu \in \Rz^2$ is given by 
	$
	g_{\bm\mu,\bm\Sigma}(\bm x)=\frac{1}{2\pi x_1x_2|\bm \Sigma|^{1/2}} \exp(-\frac{1}{2}(\log(\bm x)-\bm \mu)^T\Sigma^{-1}(\log(\bm x)-\bm \mu))\1_{\pRz^2}(\bm x)
	$  and$\log(\bm x)=(\log(x_1), \log(x_2))^T$.
	Obviously, we see that $\Mela{g_{\bm\mu,\bm\Sigma}}{\bm  c}$ is well-defined for any $\bm c\in \Rz^2$ and it holds
	\begin{align*}
	\Mela{g_{\bm\mu,\bm\Sigma}}{\bm c}(\bm t)= \exp(\bm\mu^T(\bm c-\bm 1+i\bm t))\exp\left(\frac{1}{2} (\bm c-\bm 1+i\bm t)^T \bm\Sigma (\bm c-\bm 1+i\bm t)\right), \quad \bm t\in \Rz^2.
	\end{align*}
	Then we can easily see that
	$
	|\mathcal M_{\bm c}[g_{\bm\mu, \bm\Sigma}](\bm t)|^2 = C_{\bm \mu, \bm c, \bm\Sigma} \exp(-\bm t^T\bm \Sigma\bm t)=C_{\bm \mu, \bm c, \bm\Sigma} \exp(-(t_1^2\Sigma_{11}+2t_1t_2 \Sigma_{21}+ t_2^2\Sigma_{22})), t\in \Rz.
	$
\end{example}

\begin{example}[Uniform distribution on $S:=\{(x,y)\in \pRz^2:x_1, x_2\in (1,2)\}$]\label{ex:unif}
	Consider the density $g(\bm x)= \frac{1}{3} \mathds 1_S(\bm x), \bm x\in \pRz^2$. 
	Obviously, we see that $\Mela{g}{\bm  c}$ is well-defined for any $\bm c\in \pRz^2$ and it holds
	\begin{align*}
	\Mela{g}{\bm c}(\bm t)=\frac{2^{c_1+c_2+i(t_1+t_2)}-1}{3(c_1+it_1)(c_2+it_2)}, \quad \bm t\in \Rz^2.
	\end{align*}
	Then we can see that
	$
	c_{\bm c} (1+|t_1|^2)^{-1}(1+|t_2|^2)^{-1}\leq |\mathcal M_{\bm c}[g](\bm t)|^2 \leq  C_{\bm c}  (1+|t_1|^2)^{-1}(1+|t_2|^2)^{-1}, t\in \Rz.
	$
\end{example}

The Mellin transform of the error densities in Example \ref{ex:lognorm} and \ref{ex:unif} fullfill both assumption \eqref{eq:ass:g1}. While the log normal distribution is an example for an super smooth error density, it can be seen that the density in \ref{ex:unif} is a smooth error density with $\gamma_1=\gamma_2=1$. 

\subsection{Numerical simulation}
Let us illustrate the performance of the estimator $\widehat f_{\widehat k}$ defined in \eqref{eq:est:def} and  \eqref{eq:data:driven}. We will restrict ourselves to the case $d=2$. For a simulation study with a data-driven choice of the dimension parameter for the univariate case, we refer to \cite{Brenner-MiguelComteJohannes2020} where they used a penalized contrast strategy. In the bivariate case, we will actually study the performance of the fully data-driven method presented in \eqref{eq:data:driven}, while we omit the consideration of different values of $\bm c\in \mathbb R^2$.

In the upcoming simulation study we will consider the densities
\begin{enumerate}
	\item\label{si:f:i}  \textit{Gamma Distribution}: $f_1(x)=\frac{x^{3}}{96}\exp(-0.5x)$, see Example \ref{ex:gam} and Table \ref{table:num_res2},
	\item\label{si:f:ii} \textit{Weibull Distribution}:  $f_2(x) = 2 x \exp(-x^2) \1_{\Rz^+}(x)$, see Example \ref{ex:wei} and Table \ref{table:num_res2}
	\item\label{si:f:iii}  \textit{Beta Distribution}: $f_3(x)=\frac{1}{560} (0.5x)^3(1-0.5x)^4 \1_{[0,1]}(0.5x)$, see Example \ref{ex:beta} and Table \ref{table:num_res} ,  and
	\item\label{si:f:iv} \textit{Log-Normal Distribution}: $f_4(x)=\frac{1}{\sqrt{2\pi}x} \exp(-(\log(x)^2/2))\1_{(0, \infty)}(x)$, see Example \ref{ex:lognorm} and Table \ref{table:num_res2}.
\end{enumerate}
For the error densities, we consider the univariate densities
\begin{enumerate}
	\item \label{si:err:ii} \textit{Pareto Distribution}: $g_1(x)= x^{-2} \mathds 1_{(1,\infty)}(x)$, see Tables \ref{table:num_res}, \ref{table:num_res2} and Example \ref{ex:log-gamma},
	\item\label{si:err:iv}  \textit{Log-Gamma Distribution}: $g_2(x)=\frac{1}{\Gamma(1/2)}x^{-2} \log(x)^{-1/2}\1_{(1, \infty)}(x)$, , see Tables \ref{table:num_res}, \ref{table:num_res2} and Example \ref{ex:log-gamma}.
\end{enumerate}
In the sense of \eqref{eq:ass:g1} we see that  $g_1$ has the parameter $\gamma=1$, while $g_2$ has $\gamma = 1/2$. 

Let us now consider the data-driven choice defined  \eqref{eq:est:def} and \eqref{eq:data:driven} for the case of $d=2$.  To illustrate the performance of our estimator, we consider the following three cases
\begin{enumerate}
	\item\label{si:mult:o} \textit{Error densities:} $f(x_1, x_2)= f_1(x_1)f_1(x_1)$ with direct observations compared to observations with $g(x_1,x_2)=g_2(x_1)g_2(x_2)$ with $\bm c=(1/2, 1/2)^T$,
	\item\label{si:mult:i} \textit{Anisoptopic density:}  $f(x_1, x_2)=f_3(x_1)f_4(x_2)$ with direct observations and $c=(1/2, 1/2)^T$,
	\item \label{si:mult:ii} \textit{Dependency:} $\bm X \sim \mathrm{LN}_{\mu, \bm \Sigma}$ with $\bm \mu = (\log(4), \log(4))^T$, $c=(1/2, 1/2)^T$ and direct observations. For $\bm \Sigma$ we compare
	$$\bm \Sigma_1 = \begin{pmatrix} 1& 0\\ 0 &0.81 \end{pmatrix}, \quad\Sigma_2= \begin{pmatrix} 1 &0.1 \\ 0.1 & 0.81\end{pmatrix};$$
	\item\label{si:mult:iii} \textit{Anisotropic error}: $f(\bm x)=f_2(x_1)f_2(x_2)$ with $\bm U=(U_1, U_2)^T$, where $U_1 \sim \mathrm{sL}\Gamma_{(0,1,1)}$ and $U_2\sim  \mathrm{sL}\Gamma_{(0,1/2,1)}$, dependent, and $\bm c=(1/2,1/2)^T$ 
\end{enumerate}
For the first case, we visualise the impact of observations with measurement error compared to direct observations. The second case resembles the case when the decay of the Mellin transform of the density $f$ has significantly different behaviour in different direction. The thir case shall illustrate the behaviour of the estimator when the two coordinates of $X$ are dependent while in the fourth case the decay of the Mellin transform of the density is similar but the decays of the error density are not the same.
By minimising an integrated weighted squared error over a family of histogram
densities with randomly  drawn partitions  and weights we select
$\chi= 1$ (respectively 
$\chi = 0.3$) for the cases of direct observation (respectively contaminated data), where $\chi$ is the variance constant, see \eqref{eq:data:driven}.
\paragraph{Case 1:}
Let us start by compare the influence of measurement errors by comparing the estimator $\widehat f_{\widehat k}$ based on the copies of $\bm X$ compared to copies of $\bm Y$.

\begin{minipage}{\textwidth}
	\begin{minipage}[t]{0.49\textwidth}
		\includegraphics[width=\textwidth,height=35mm]{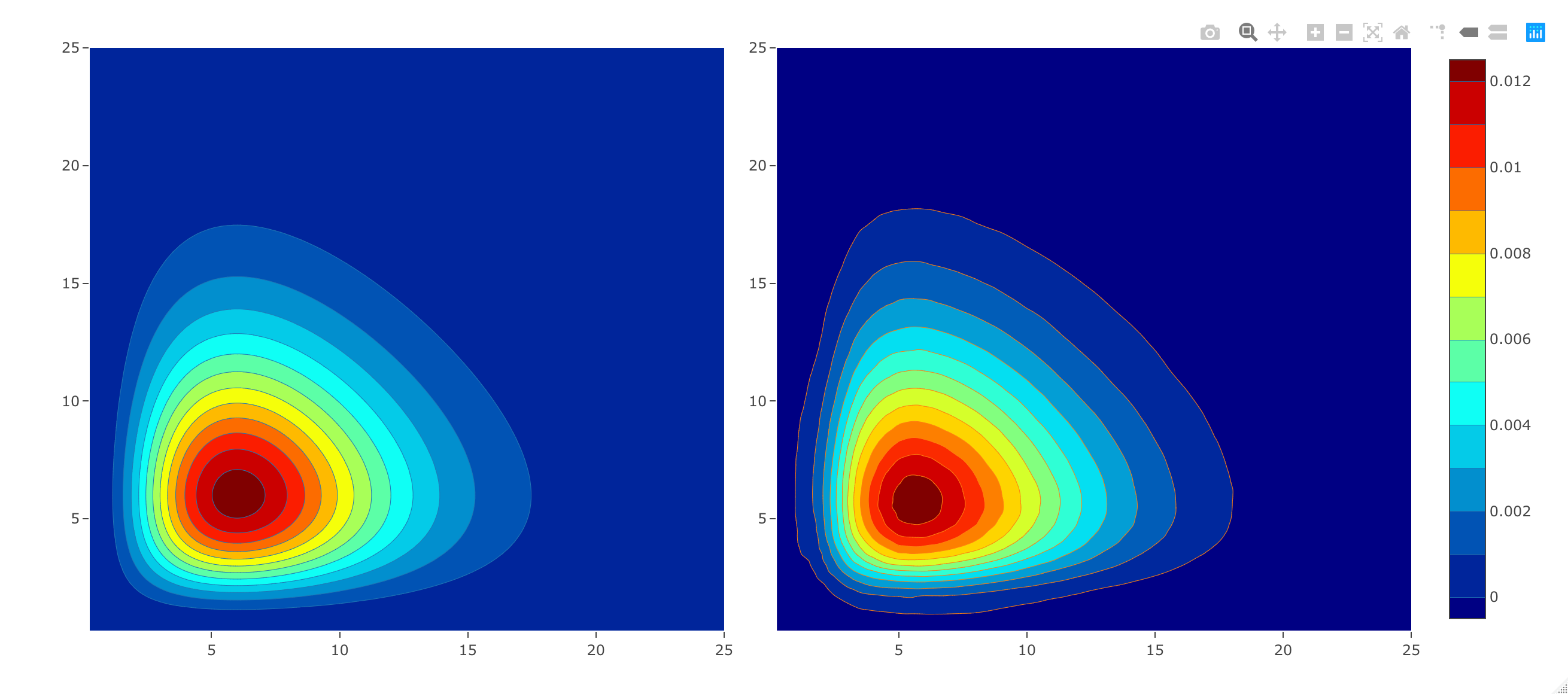}
	\end{minipage}
	\begin{minipage}[t]{0.49\textwidth}
		\includegraphics[width=\textwidth,height=35mm]{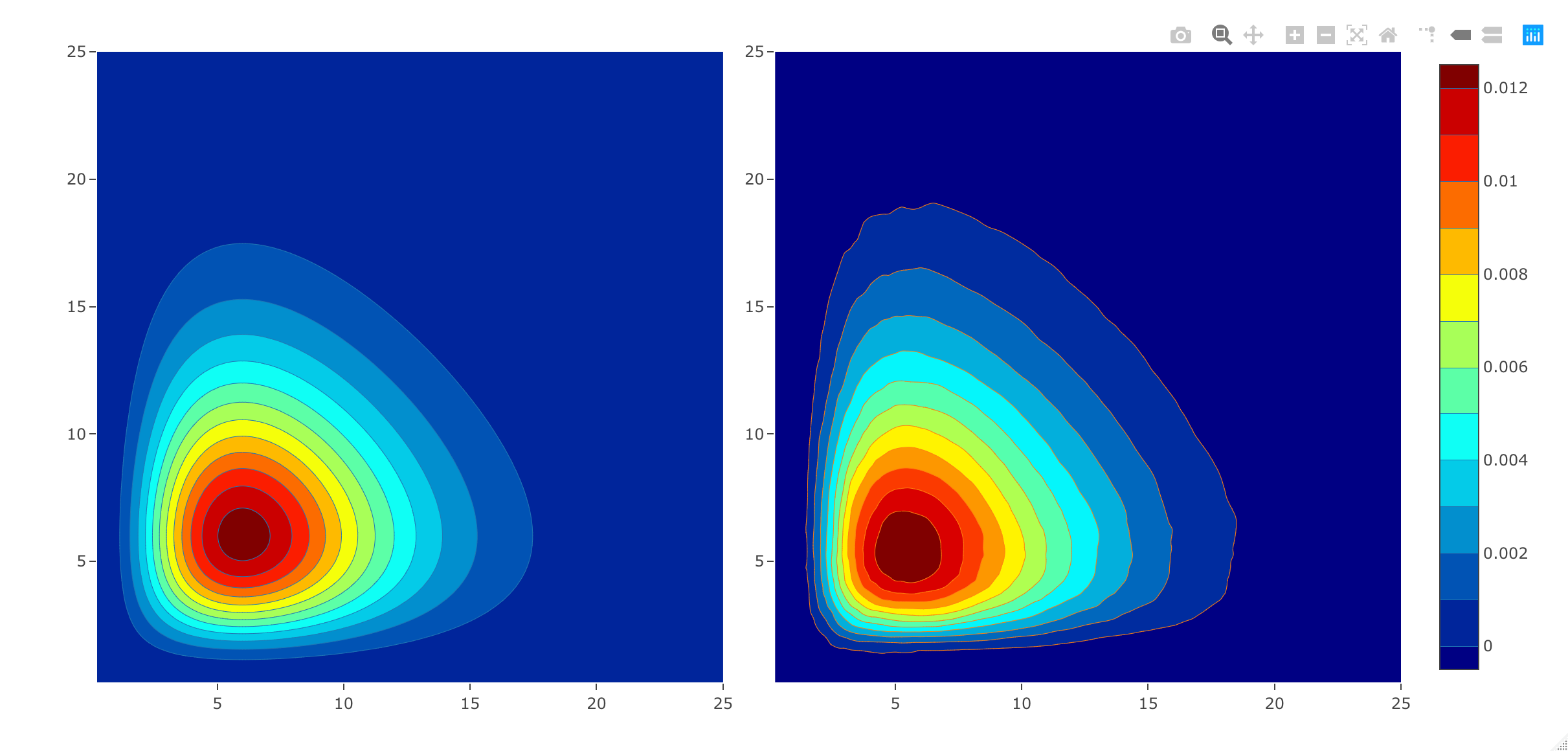}
	\end{minipage}
	\begin{minipage}[t]{\textwidth}
		\begin{minipage}[t]{0.24\textwidth}
			\includegraphics[width=\textwidth,height=35mm]{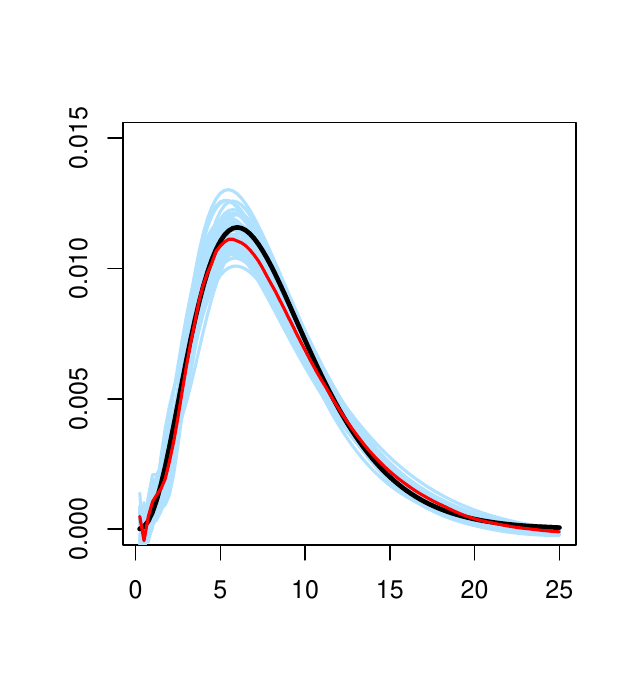}
		\end{minipage}
		\begin{minipage}[t]{0.24\textwidth}
			\includegraphics[width=\textwidth,height=35mm]{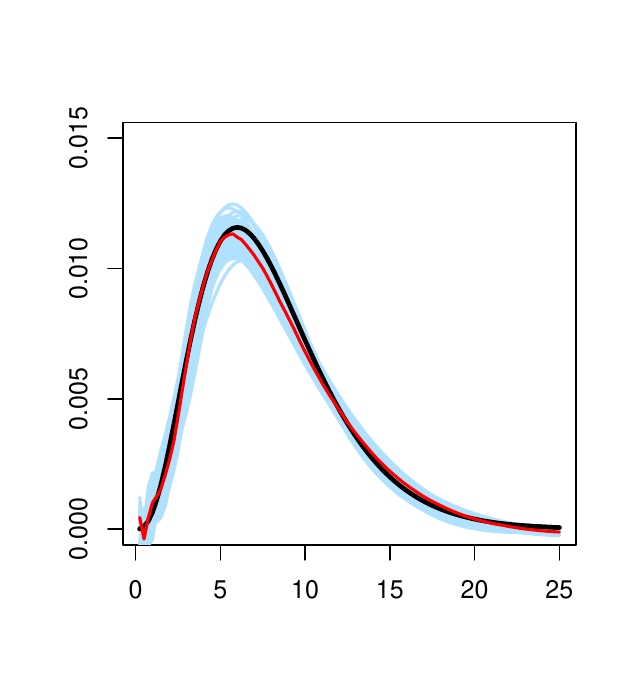}
		\end{minipage}
		\begin{minipage}[t]{0.24\textwidth}
			\includegraphics[width=\textwidth,height=35mm]{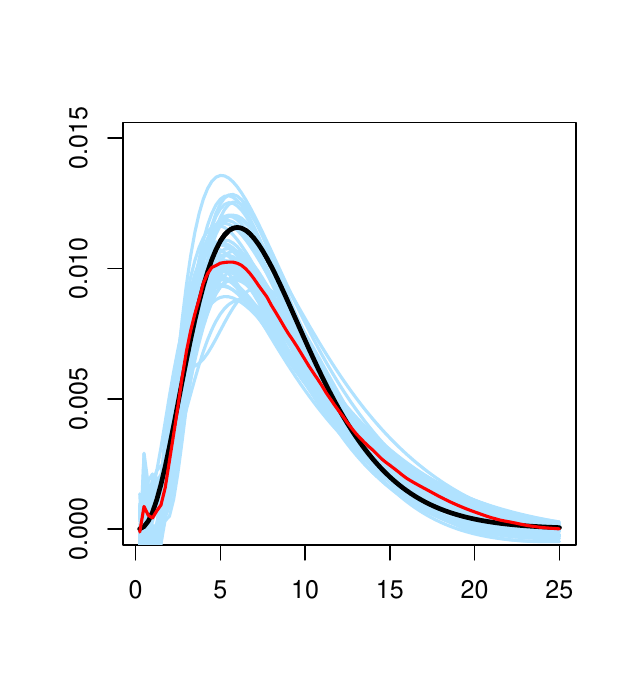}
		\end{minipage}
		\begin{minipage}[t]{0.24\textwidth}
			\includegraphics[width=\textwidth,height=35mm]{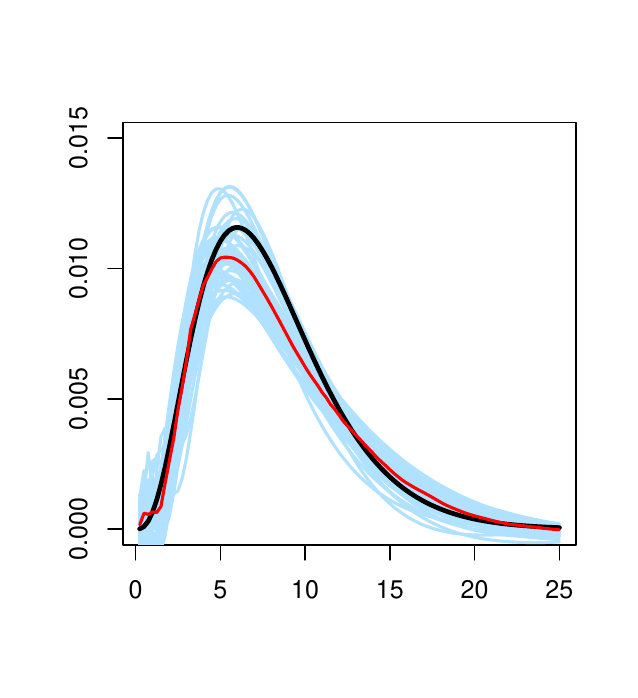}
		\end{minipage}
	\end{minipage}
	\captionof{figure}{\label{figure:1}The estimator $\widehat f_{\widehat {\bm k}}$ are depicted for 
		50  Monte-Carlo simulations with  sample size $n=1000$ with direct observations (left) and with multiplicative measurement errors (right). The top plots are the true density (left) and the pointwise median of the estimators (right). The bottom plots are the sections for $x=7.5$ (right) and $y=7.5$ (left) where the true density $f$ is given by the black curve while the red curve is the pointwise empirical median of the 50 estimates.}
\end{minipage}\\[2ex]

\paragraph{Case 2:}
In the second case, we additionally compare the anisotropic estimator $\widehat f_{\widehat{\bm k}}$ with the isotropic choice, that is we define $\widehat f_{\widetilde k}:=\widehat f_{\widetilde k, \widetilde k}$ with 
\begin{align} \widetilde k:= \argmin_{k\in \mathbb N, \Delta_g((k,k)) \leq n} -\| \widehat f_{(k,k)}\|_{x^{2c-1}}^2 + \kappa \widehat{\sigma} \frac{\Delta_g((k,k))}{n}, \label{eq:iso}
\end{align}
a penalized contrast upproach which is a direct generalisation of the estimator given in \cite{Brenner-MiguelComteJohannes2020}. Here we choose $\kappa=5$ by a preliminary simulation study.\\
\begin{minipage}{\textwidth}
	\begin{minipage}[t]{0.49\textwidth}
		\includegraphics[width=\textwidth,height=35mm]{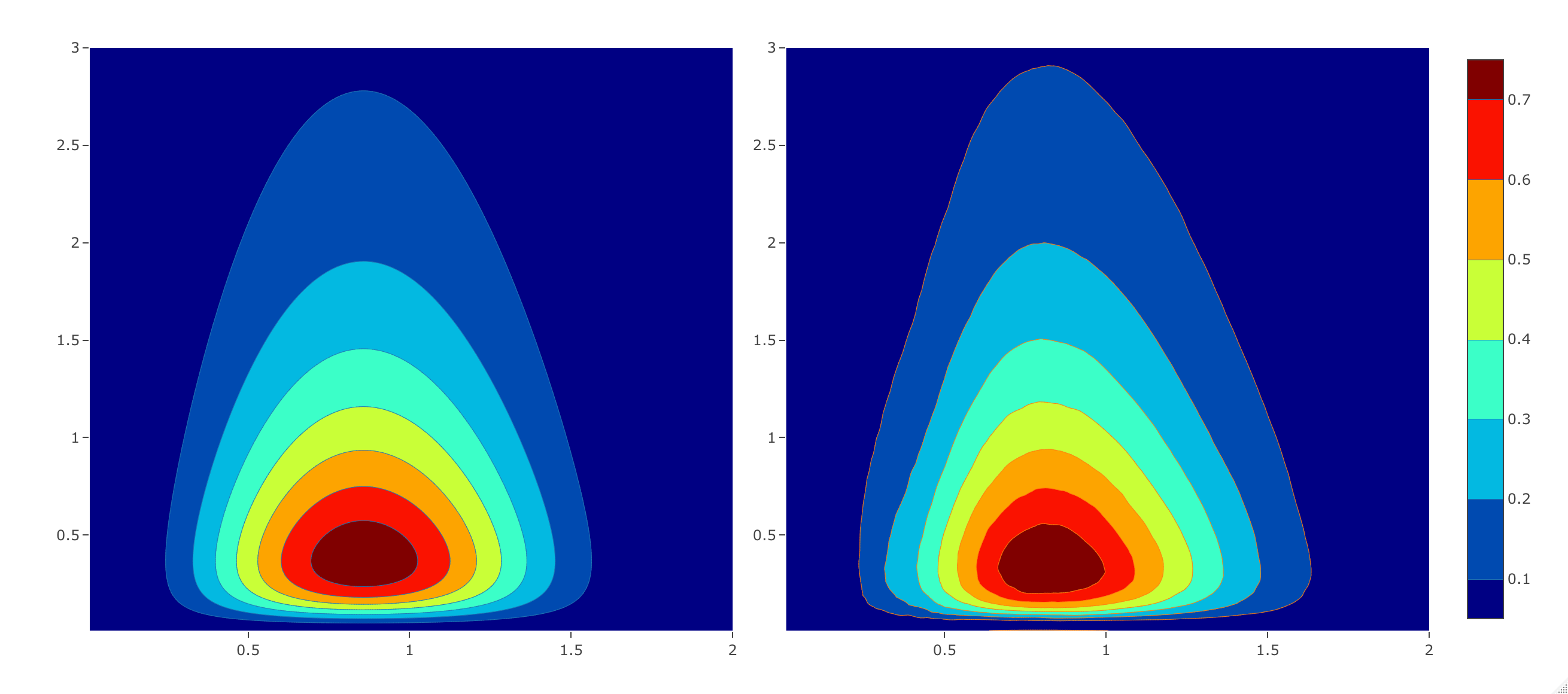}
	\end{minipage}
\begin{minipage}[t]{0.49\textwidth}
	\includegraphics[width=\textwidth,height=35mm]{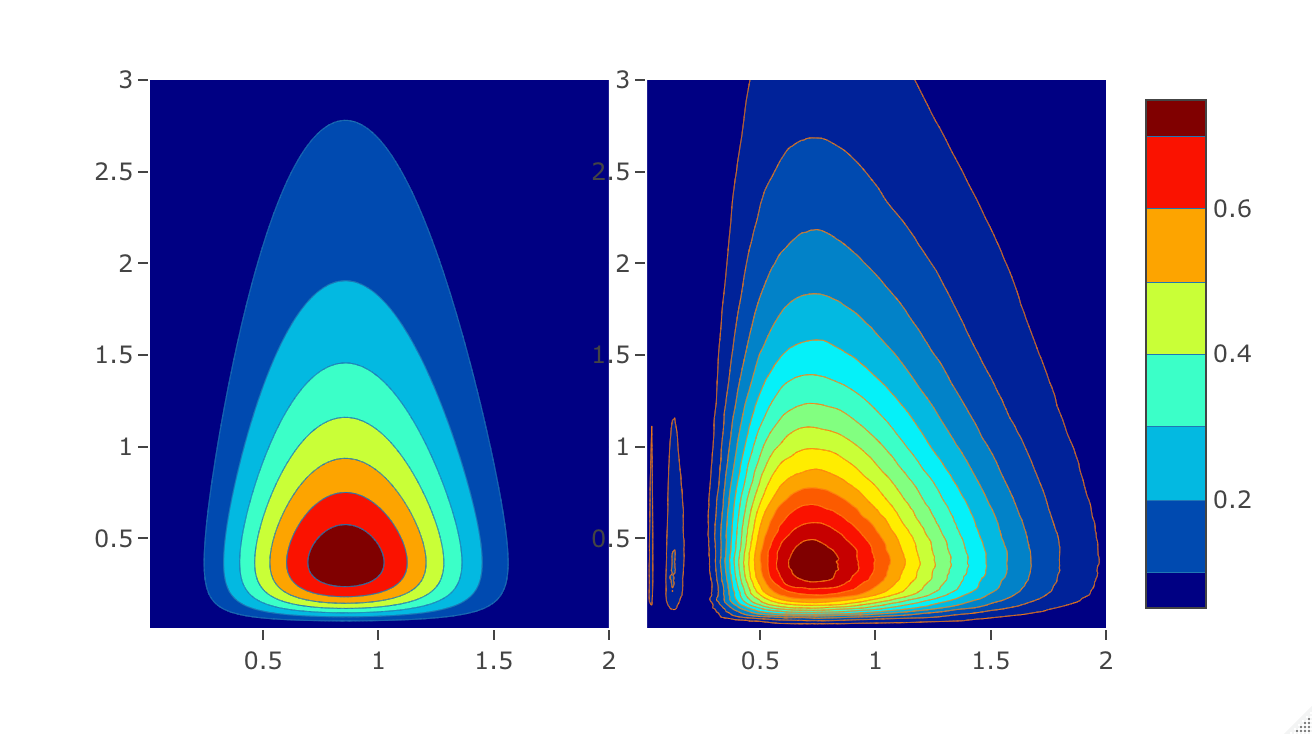}
\end{minipage}
	\begin{minipage}[t]{\textwidth}
		\begin{minipage}[t]{0.24\textwidth}
		\includegraphics[width=\textwidth,height=35mm]{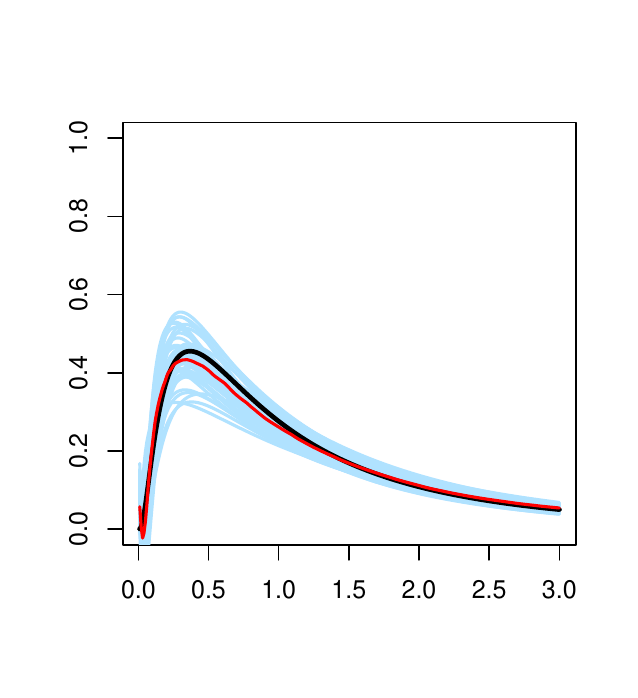}
	\end{minipage}
\begin{minipage}[t]{0.24\textwidth}
		\includegraphics[width=\textwidth,height=35mm]{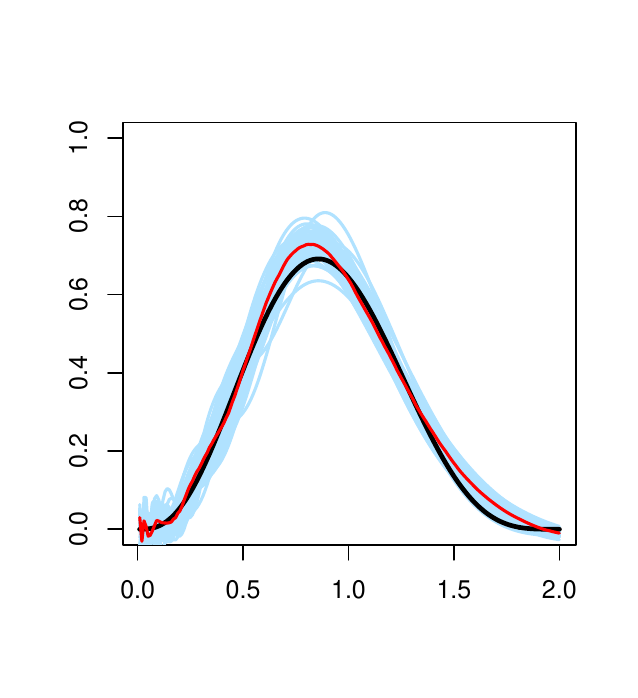}
	\end{minipage}
\begin{minipage}[t]{0.24\textwidth}
	\includegraphics[width=\textwidth,height=35mm]{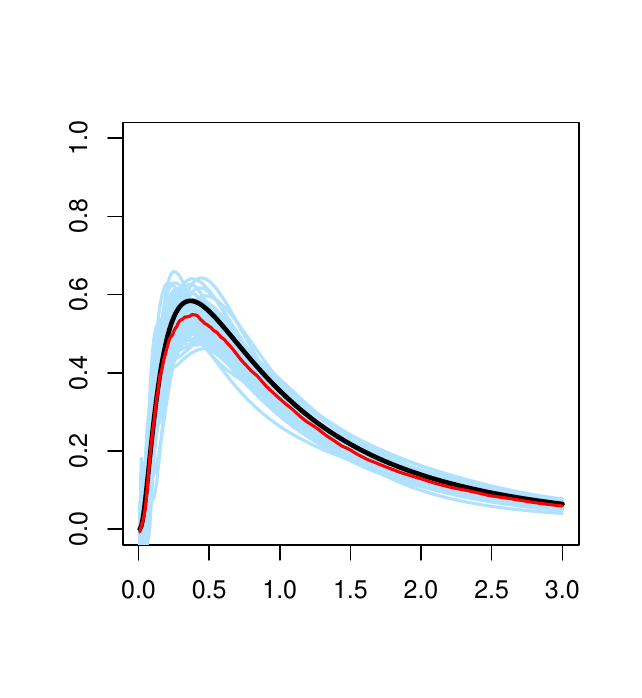}
\end{minipage}
\begin{minipage}[t]{0.24\textwidth}
	\includegraphics[width=\textwidth,height=35mm]{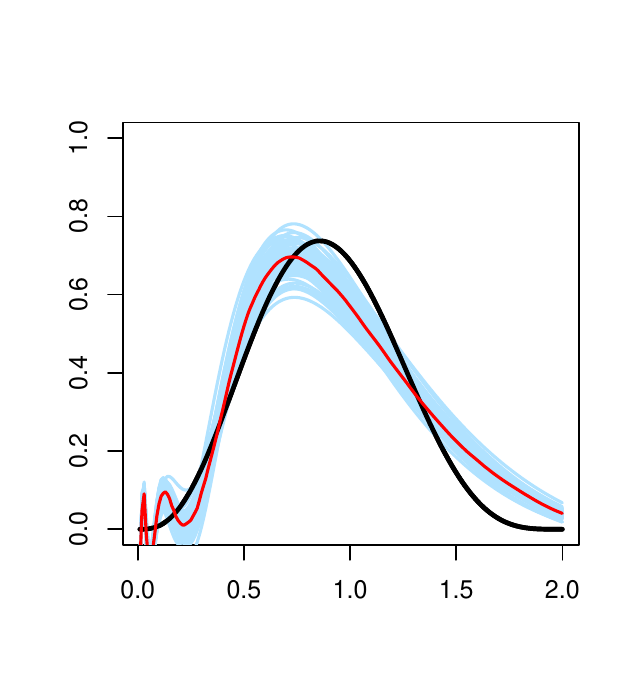}
\end{minipage}
\end{minipage}
	\captionof{figure}{\label{figure:5}The estimator $\widehat f_{\widehat {\bm k}}$ (left) and $\widehat f_{\widetilde k}$ (right) are depicted for 
		50  Monte-Carlo simulations with  sample size $n=1000$ with direct observations. The top plots are the true density (left) and the pointwise median of the estimators (right). The bottom plots are the sections for $x=5$ (right) and $y=0.59$ (left) where the true density $f$ is given by the black curve while the red curve is the pointwise empirical median of the 50 estimates.}
\end{minipage}\\[2ex]
As one can see in Fig. \ref{figure:5}, the anisotropic estimator seems to invest more in the approximation of the Beta distribution than the Log normal distribution. This leads to worse performance in the Log normal direction but to an overall satisfying result. In comparison to that, the isotropic estimator chooses in both direction the same cut-off parameter leading to a better approximation of the log normal distribution but also to a worse approximation of the beta distribution. Overall it seems that the anisotropic estimator behaves better.
\paragraph{Case 3:} 
Now we consider the influence of the dependency between  the coordinates of $\bm X$. While $\bm \Sigma_1$ resembles the case of independent coordinates, $\bm\Sigma_2$ is not a diagonal matrix and thus the coordinates are dependent.\\
\begin{minipage}{\textwidth}
	\begin{minipage}[t]{0.49\textwidth}
		\includegraphics[width=\textwidth,height=35mm]{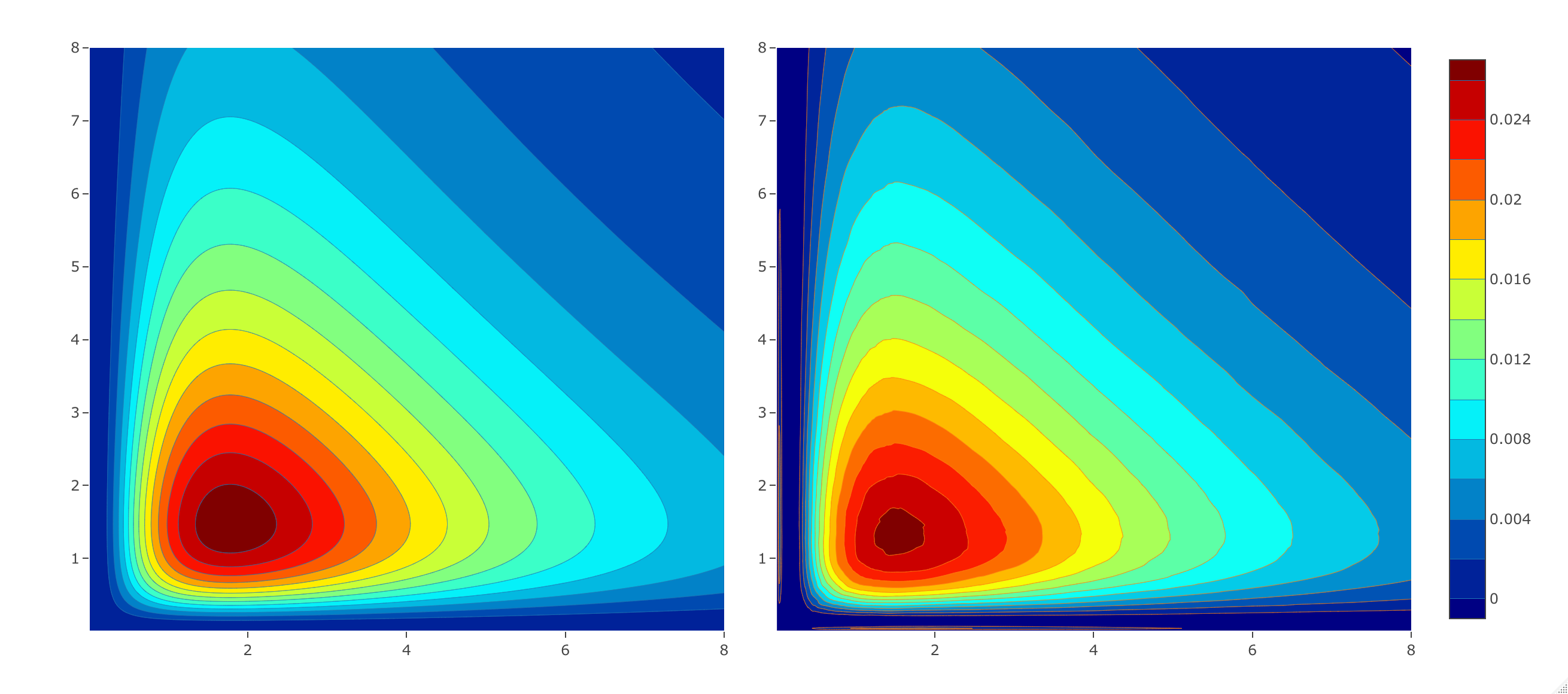}
	\end{minipage}
	\begin{minipage}[t]{0.49\textwidth}
		\includegraphics[width=\textwidth,height=35mm]{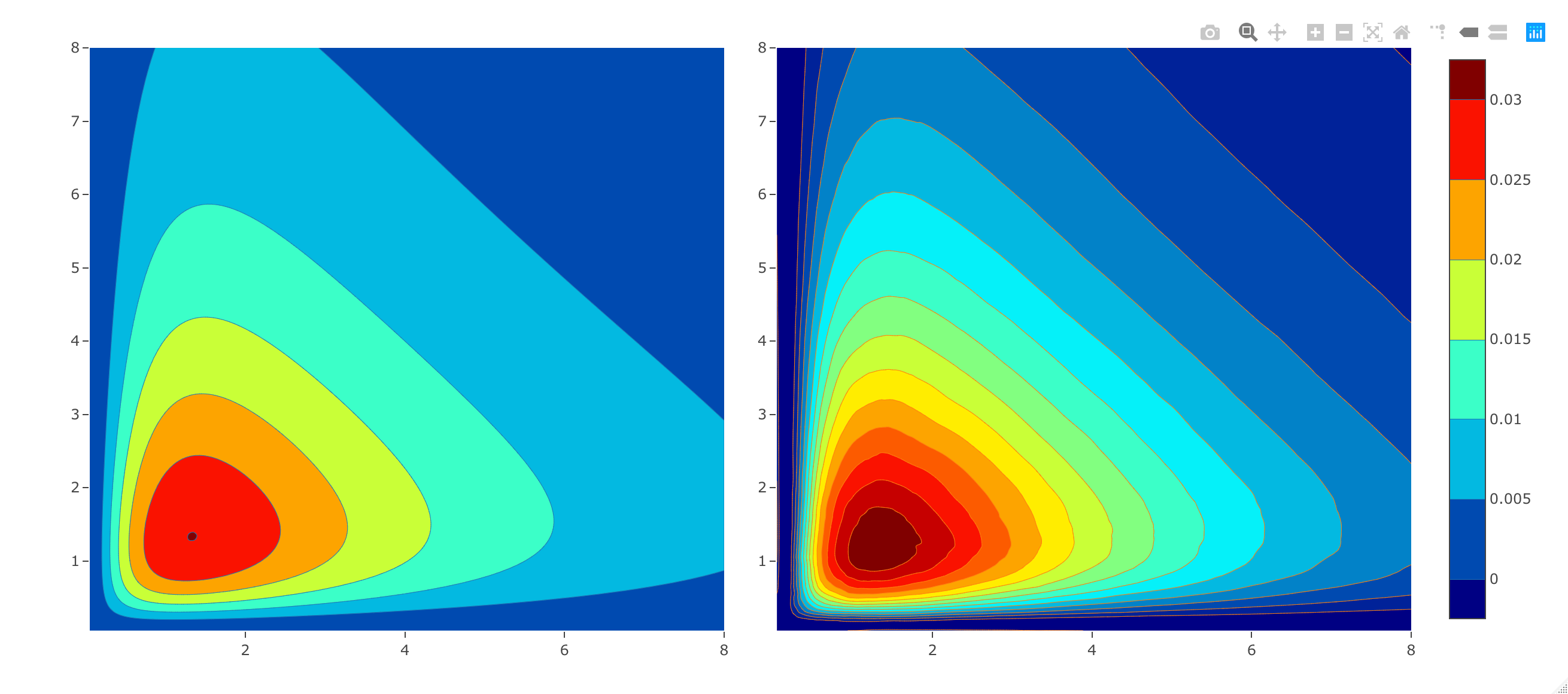}
	\end{minipage}
	\begin{minipage}[t]{\textwidth}
		\begin{minipage}[t]{0.24\textwidth}
			\includegraphics[width=\textwidth,height=35mm]{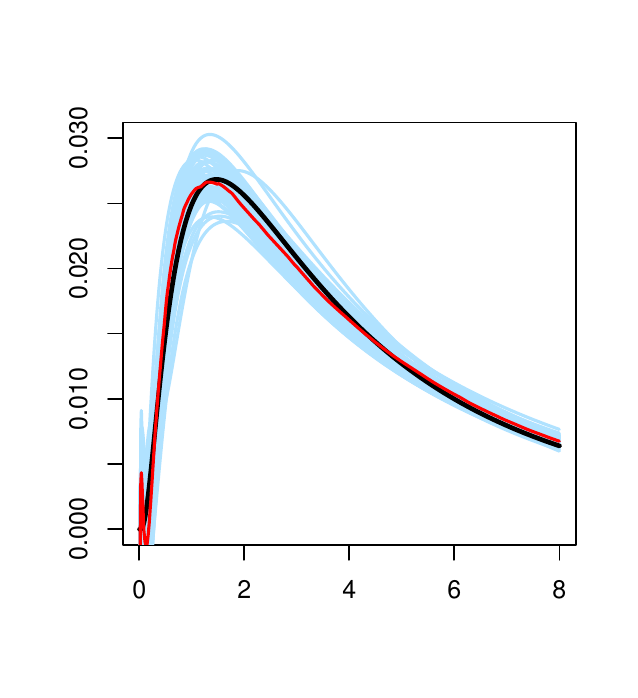}
		\end{minipage}
		\begin{minipage}[t]{0.24\textwidth}
			\includegraphics[width=\textwidth,height=35mm]{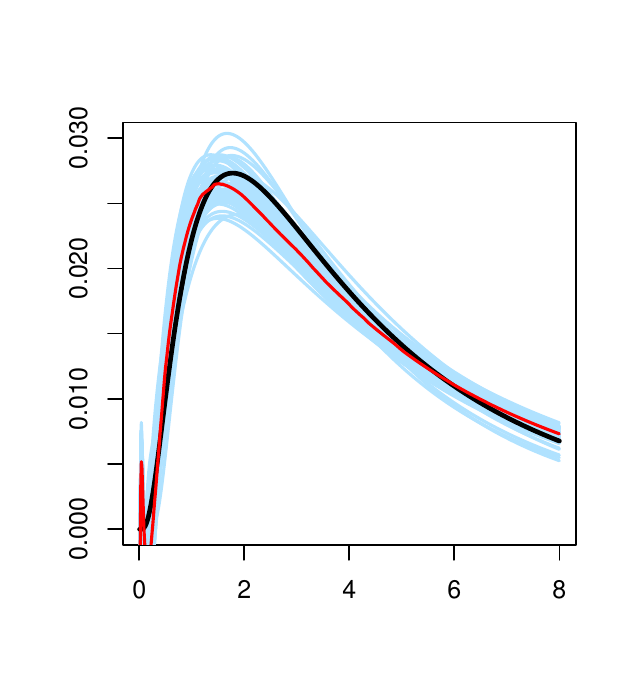}
		\end{minipage}
		\begin{minipage}[t]{0.24\textwidth}
			\includegraphics[width=\textwidth,height=35mm]{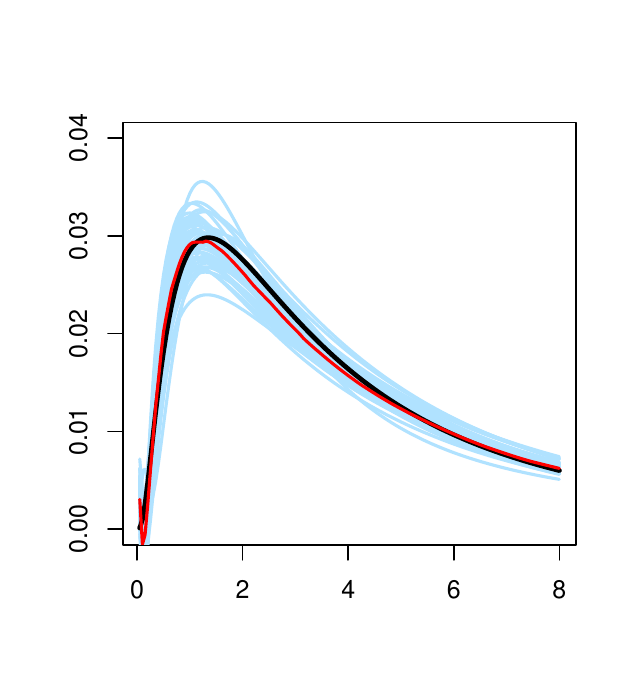}
		\end{minipage}
		\begin{minipage}[t]{0.24\textwidth}
			\includegraphics[width=\textwidth,height=35mm]{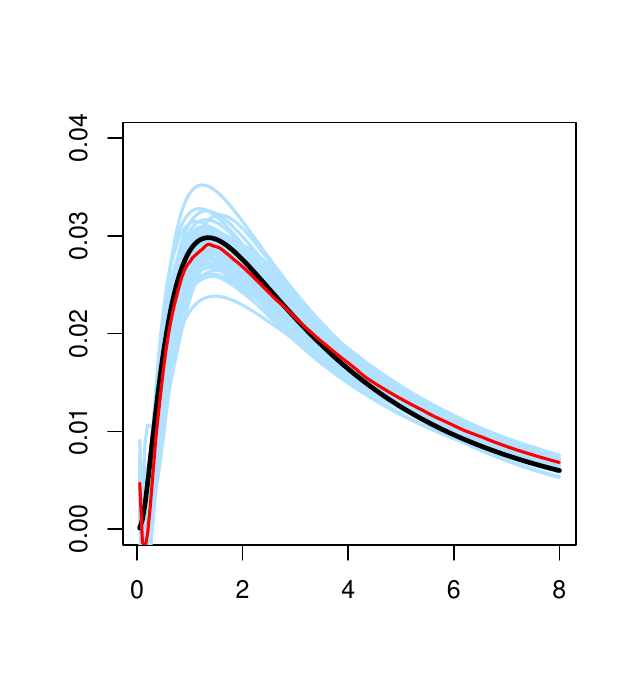}
		\end{minipage}
	\end{minipage}
	\captionof{figure}{\label{figure:6}The estimator $\widehat f_{\widehat {\bm k}}$  are depicted for 
		50  Monte-Carlo simulations with  sample size $n=1000$ with direct observations and $\bm\Sigma_1$ (left) and $\bm \Sigma_2$ (right). The top plots are the true density (left) and the pointwise median of the estimators (right). The bottom plots are the sections for $x=1.5$ (right) and $y=1.5$ (left) where the true density $f$ is given by the black curve while the red curve is the pointwise empirical median of the 50 estimates.}
\end{minipage}\\[2ex]

In Fig. \ref{figure:6} we can see that allthough the estimator does reconstruct the general shape of the density $f$, the included dependency impede slightly the estimation.

\paragraph{Case 4:}
We finish our simulation study by considering the case, where the decay of the Mellin transform of the density behaves similar in both direction, while the decay of the Mellin transform of the error densities differs. For $\bm U=(U_1, U_2)^T$ we set $U_1:=\xi_1\xi_2$ and $U_2:=\xi_2$ where $\xi_1, \xi_2 \stackrel{i.i.d.}{\sim} \mathrm{sL\Gamma}_{(0, 1/2, 1)}$, compare Example \ref{ex:log-gamma}. Then we have $U_1 \sim \mathrm{sL}\Gamma_{(0, 1, 1)}$ and $U_2 \sim \mathrm{sL}\Gamma_{(0, 1/2, 1)}$ and $$ \mathcal M_{\bm{1/2}}[g](\bm t)= (3/2 -it_1)^{-1/2} (3/2-i(t_1+t_2))^{-1/2},\quad  \bm t\in \Rz^2,$$ leading that $g$ satisfies in this situation \eqref{eq:ass:g1} with $\bm \gamma = (1, 1/2)$. For this case, we deduced by a preliminary simulation study the choice $\chi_1=\chi_2=0.25$. For the distribution of $\bm X=(X_1, X_2)^t$ we set for the sake of simplicity $X_1, X_2 \stackrel{i.i.d.}{\sim} W_{(2,1)}$, see Example \ref{ex:wei}. We compare the performance of the data-driven anisotropic estimator $\widehat f_{\widehat k}$ with the performance of the data-driven isotropic estimator $\widehat f_{\widetilde k}$ introduced in \eqref{eq:iso} with the choice $\kappa = 1.02$. For both estimator we consider $\bm c=\bm{1/2}=(1/2, 1/2)^T.$

\begin{minipage}{\textwidth}
	\begin{minipage}[t]{0.49\textwidth}
		\includegraphics[width=\textwidth,height=35mm]{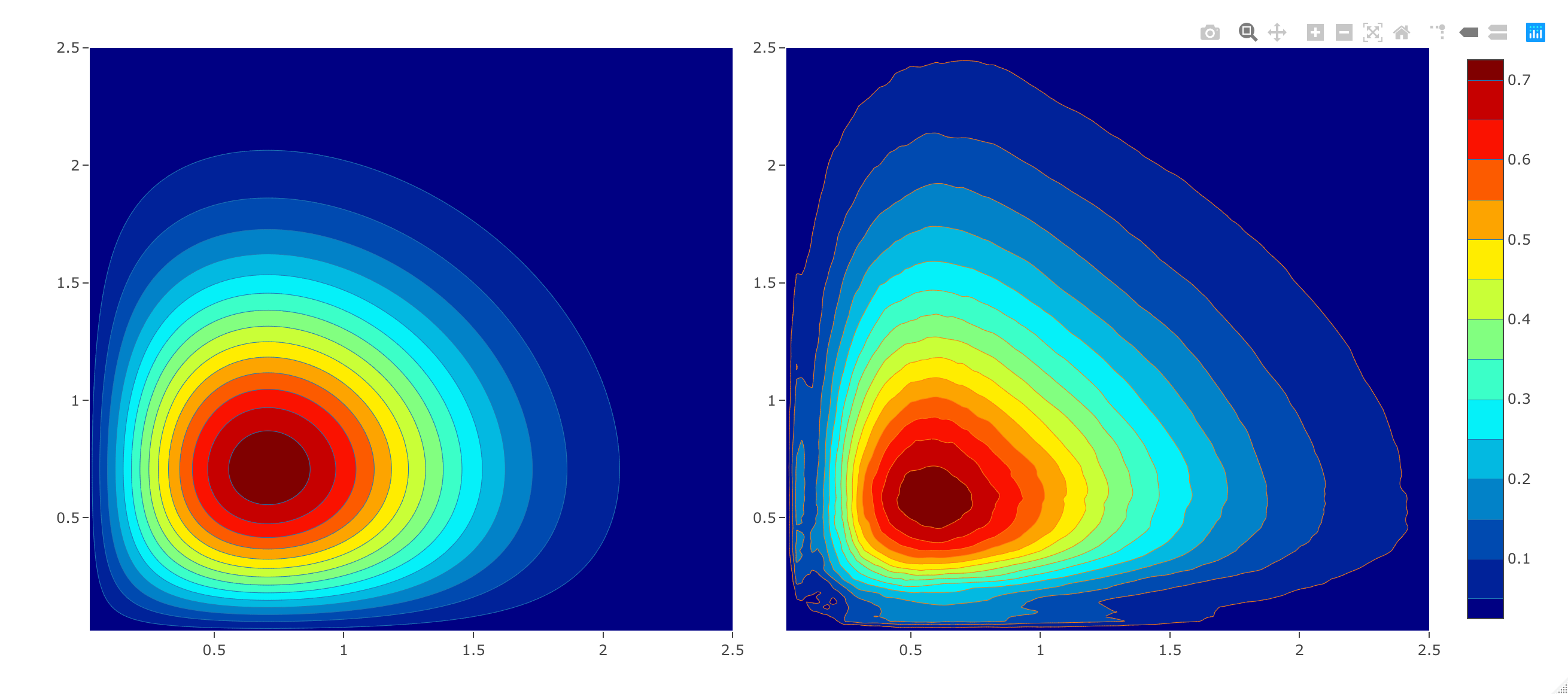}
	\end{minipage}
	\begin{minipage}[t]{0.49\textwidth}
		\includegraphics[width=\textwidth,height=35mm]{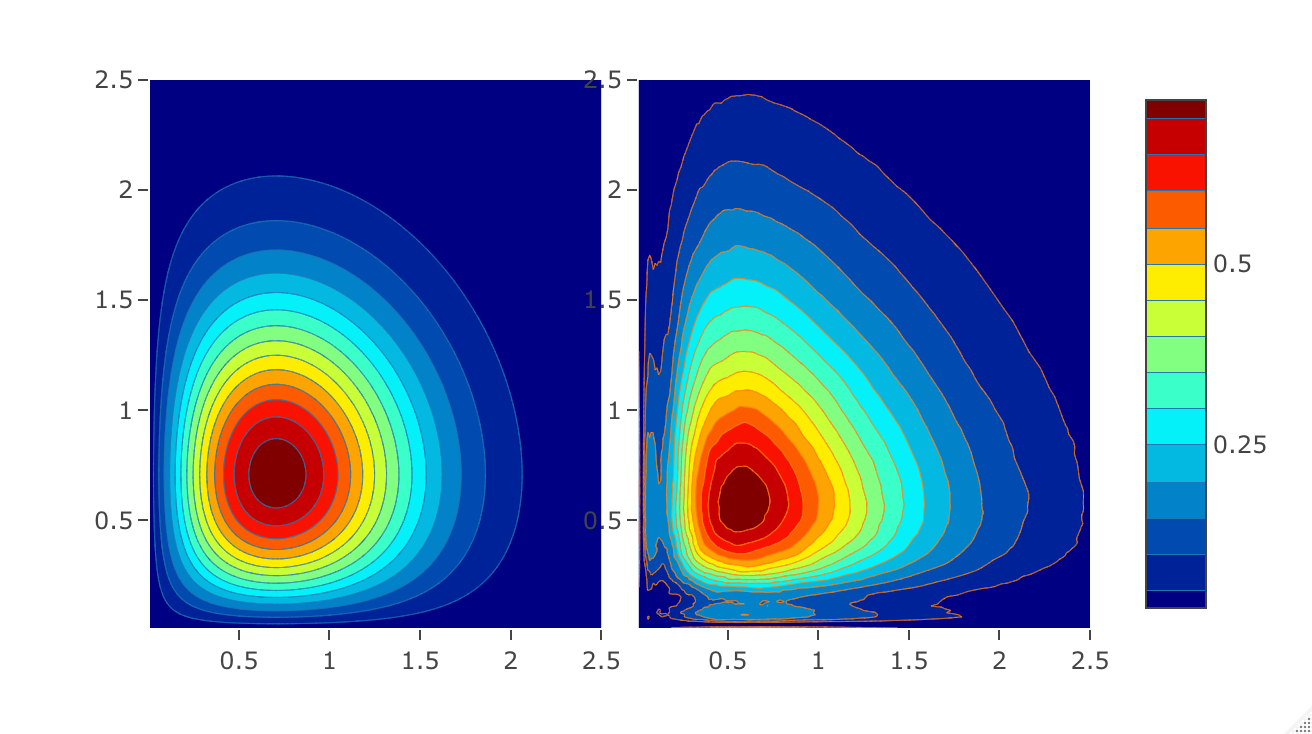}
	\end{minipage}
	\begin{minipage}[t]{\textwidth}
		\begin{minipage}[t]{0.24\textwidth}
			\includegraphics[width=\textwidth,height=35mm]{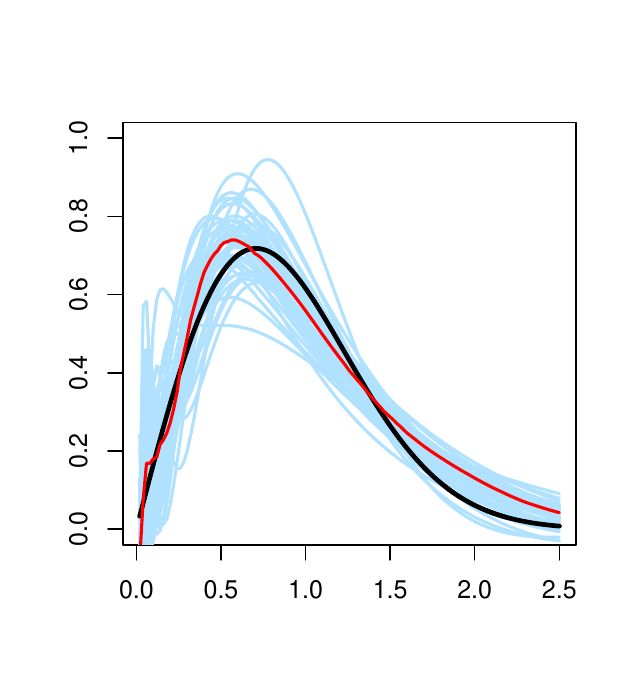}
		\end{minipage}
		\begin{minipage}[t]{0.24\textwidth}
			\includegraphics[width=\textwidth,height=35mm]{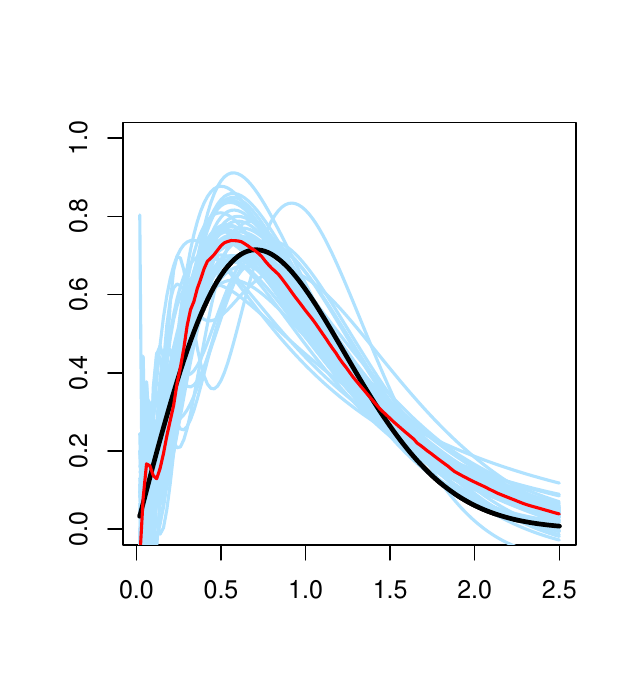}
		\end{minipage}
		\begin{minipage}[t]{0.24\textwidth}
			\includegraphics[width=\textwidth,height=35mm]{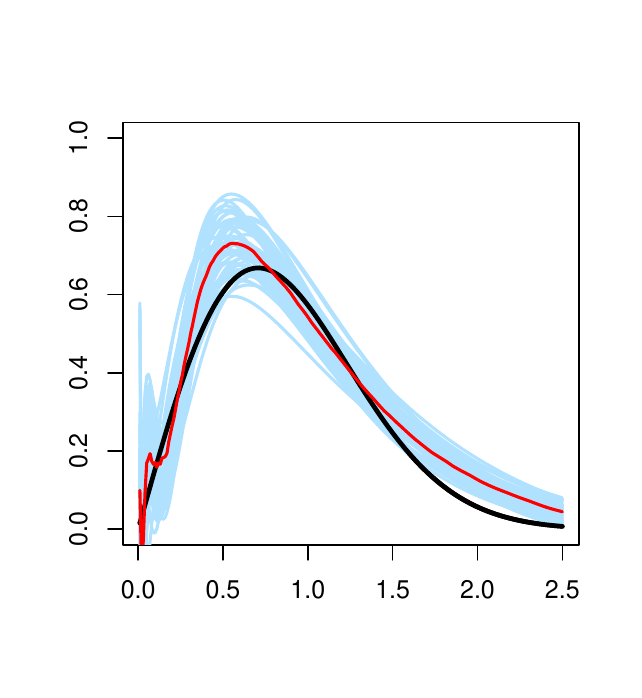}
		\end{minipage}
		\begin{minipage}[t]{0.24\textwidth}
			\includegraphics[width=\textwidth,height=35mm]{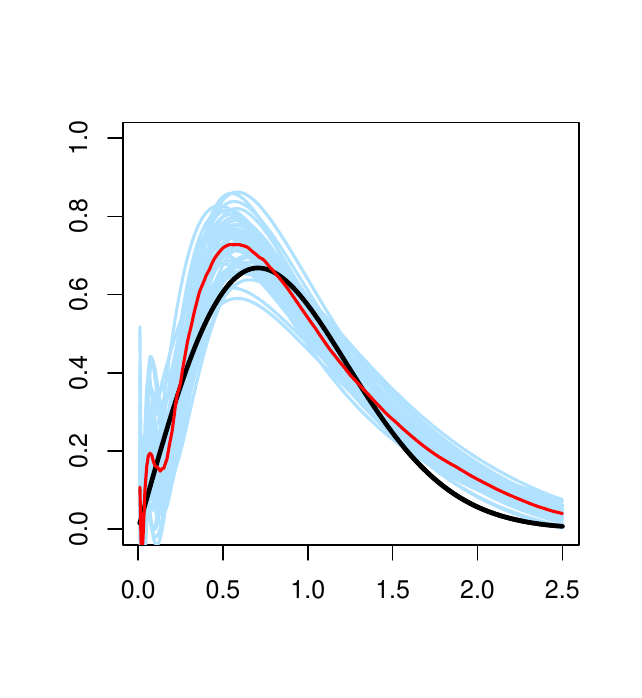}
		\end{minipage}
	\end{minipage}
	\captionof{figure}{\label{figure:7}The estimator $\widehat f_{\widehat {\bm k}}$ (left) and $\widehat f_{\widetilde k}$ (right) are depicted for 
		50  Monte-Carlo simulations with  sample size $n=1000$ with observations under multiplicative measurement errors. The top plots are the true density (left) and the pointwise median of the estimators (right). The bottom plots are the sections for $x=5$ (right) and $y=0.59$ (left) where the true density $f$ is given by the black curve while the red curve is the pointwise empirical median of the 50 estimates.}
\end{minipage}\\[2ex]

As one can see in Fig. \ref{figure:7} the anisotropic estimator $\widehat f_{\widetilde k}$ behaves better in the second coordinates as the isotropic estimator $\widehat f_{\widetilde k}$ , which is consistent with the theory since the decay of the Mellin transform of the error density in this direction is slower. 

\subsection{Comment}

As seen by the simulation study, the anisotropic data-driven estimator $\widehat f_{\widehat k}$ behaves reasonable for the case $d=2$. Compared to the an isotropic estimator, the desired flexibility of an anisotropic choice can be interpretated in the Figures \ref{figure:5} and \ref{figure:7}, where an anisotropic choice is necessary due to different behaviour of the decay of the Mellin transform of the density $f$, respectively the error density $g$. Allthough the implementation of the estimator for the case $d>2$ is possible, we do not provide a simulation study for these cases. It is worth stressing out, that for higher dimensions, the estimator does suffers under the well-known curse of dimensionality, that is, the convergence rates, and therefore the performance of the estimation strategy, is slower than in lower dimensions. 
\newpage

\section*{Acknowledgments}

I thank the Editor, Associate Editor and referees for their helpful comments and valuable suggestions. This research was supported by the Research Training Group ”Statistical Modeling of Complex Systems” funded by the German Science Foundation. Furthermore, I want to thank Jan Johannes and Fabienne Comte  for their helpful advice and their support of my work. 

\section*{Appendix}

\subsection{Preliminaries}\label{a:prel}
%

% --------------------------------------------------------------------
% <<Re \ref{re:diff_fun}>>
% --------------------------------------------------------------------

We will now present some proof sketches for the properties of the Mellin transform stated in Section \ref{mt}. We recall that for a function $h\in \Lz^2(\pRz^d, \bm x^{\ushort{2\bm c-\bm 1}})$ we defined the Mellin transform developed in $c\in \Rz$ by
\begin{align*}
\Mela{h}{\bm c}(\bm t):= (2\pi)^{d/2} \mathcal F[\Phi_{\bm c}[h]](\bm t), \quad \bm t\in \Rz^d,
\end{align*} 
where 
\begin{align*}
\hspace*{-1cm}\mathcal F: \Lz^2(\Rz^d)\rightarrow \Lz^2(\Rz^d), H\mapsto (\bm t\mapsto \mathcal F[H](\bm t):=\lim_{k\rightarrow \infty}(2\pi)^{-d/2}\int_{B_k(\bm 0)} \exp(- i \langle\bm t, \bm x\rangle) H(\bm x) d\bm t)
\end{align*}
is the Plancherel-Fourier transform where the limit is understood in a $\Lz^2(\Rz^d)$ convergence sense and   the function $\Phi_{\bm c}$ is defined by $\Phi_{\bm c}:\Lz^2(\pRz^d, \bm x^{\ushort{2\bm c-\bm 1}})\rightarrow \Lz^2(\Rz^d), h\mapsto \varphi^{\ushort{\bm c}} h\circ \varphi$ and $\varphi:\Rz^d \rightarrow \pRz^d$,  $\bm x\mapsto (\exp(- x_1), \dots, \exp(- x_d))^T$.

By assuming $h\in \Lz^2(\pRz^d, \bm x^{\ushort{2\bm c- \bm1}})\cap \Lz^1(\pRz^d,\bm x^{\ushort{\bm c-\bm1}})$ we get that $\Phi_{\bm c}[h] \in \Lz^1(\Rz^d) \cap \Lz^2(\Rz^d)$. In that case, we know that we can write the Fourier-Plancherel transform explicitly and get for all $\bm t\in \Rz^d$,
\begin{align*}
\Mela{h}{\bm c}(\bm t)&= \int_{\Rz^d}\exp(- i \langle \bm t, \bm x\rangle ) \varphi(\bm x)^{\ushort{\bm c}}h(\varphi(\bm x)) d\bm x= \int_{\Rz^d}\varphi(\bm x)^{\ushort{i\bm t}}\varphi(\bm x)^{\ushort{\bm c}}h(\varphi(\bm x)) d\bm x
= \int_{\pRz^d} \bm x^{\ushort{\bm c-\bm 1+i\bm t}} h(\bm x) d\bm x
\end{align*}
by single change of variables. 

Since $\mathcal M_{\bm c}: \Lz^2(\pRz^d,\bm x^{\ushort{2\bm c-\bm1}}) \rightarrow \Lz^2(\Rz^d)$ is a composition of isomorphism we see that it is invertible and its inverse $\mathcal M_{\bm c}^{-1}: \Lz^2(\Rz^d)\rightarrow \Lz^2(\pRz^d,\bm x^{\ushort{2\bm c-\bm 1}})$ can be expressed through $\mathcal M_{\bm c}^{-1}[H]= (2\pi)^{-d/2} \Phi_{\bm c}^{-1}[\mathcal F^{-1}[H]]$ for any $H\in \Lz^2(\Rz^d)$ where $\mathcal F^{-1}$ is the inverse of the Fourier-Plancherel transform. If additionally $H\in \Lz^1(\Rz^d)$, we can express the inverse Fourier-Plancherel transform explicitly and get for any $x\in \pRz$,
\begin{align*}
\mathcal M_{\bm c}^{-1}[H](\bm x)=(2\pi)^{-d} \Phi_{\bm c}^{-1}[\int_{\Rz^d} \exp( i \langle\bm t, \bm x\rangle) H(\bm t) d\bm t ] = \frac{1}{(2\pi)^d} \int_{\Rz^d} \bm x^{\ushort{-\bm c-i\bm t}} H(\bm t) d\bm t.
\end{align*}
Next, we are going to show a Plancherel-type equation for the Mellin transform, that is, for any $h_1,h_2 \in \Lz^2(\pRz^d, \bm x^{\ushort{2\bm c-\bm 1}})$ holds $\langle h_1, h_2\rangle_{\bm x^{\ushort{2\bm c-\bm 1}}} = (2\pi)^{-d} \langle \Mela{h_1}{\bm c}, \Mela{h_2}{\bm c} \rangle_{\Rz^d}$. Again we see that by application of a change of variable and the Plancherel equation for the Fourier-Plancherel transform that
\begin{align*}
\langle h_1, h_2 \rangle_{\bm x^{\ushort{2\bm c-\bm 1}}} &= \int_{\pRz^d} \bm x^{\ushort{\bm c}} h_1(\bm x) \bm x^{\ushort{\bm c}} h_2(\bm x) \bm x^{\ushort{-\bm 1}} d\bm x 
=\langle \Phi_{\bm c}[h_1], \Phi_{\bm c}[h_2] \rangle_{\Rz^d} \\
&=\langle \mathcal F[\Phi_{\bm c}[h_1]], \mathcal F[\Phi_{\bm c}[h_2]] \rangle_{\Rz^d} = (2\pi)^{-d} \langle \Mela{h_1}{\bm c}, \Mela{h_2}{\bm c}\rangle_{\Rz^d}.
\end{align*}
Now let us finish this short introduction by showing the convolution theorem for the Mellin transform, that is for $h_1, h_2 \in \Lz^1(\pRz^d,\bm x^{\ushort{\bm c-\bm 1}}) $ holds $\mathcal M_{\bm c}[(h_1*h_2)]=\Mela{h_1}{\bm c}\Mela{h_2}{\bm c}$ where $(h_1*h_2)$ denotes the multiplicative convolution of $h_1$ and $h_2$ which was given by
\begin{align*}
(h_1*h_2)(\bm y)= \int_{\pRz^d} h_1(\bm x) h_2(\bm y/\bm x)\bm x^{\ushort{\minus\bm 1}}d\bm x.
\end{align*}
We first show that this expression is well-defined for any $h_1,h_2 \in \Lz^1(\pRz^d, \bm x^{\ushort{\bm c-\bm1}})$.
Since  $h_1, h_2 \in \Lz^1(\pRz^d, \bm x^{\ushort{\bm c-\bm1}}) $ we have $\Phi_{\bm c}[h_1], \Phi_{\bm c}[h_2]\in \Lz^1(\Rz^d) $. Defining for two functions $H_1, H_2\in \Lz^1(\Rz^d) \cap \Lz^2(\Rz^d)$ the additive convolution $(H_1*_+H_2)$ by $(H_1*_+ H_2)(\bm y):= \int_{\Rz^d} H_1(\bm y-\bm x) H_2(\bm x) d\bm x$ for any $\bm y\in \Rz^d$ we know from functional analysis that $(H_1 *_+ H_2) \in \Lz^1(\Rz^d) $. Thus it follows that $(h_1*h_2)=\Phi_{\bm c}^{-1}[(\Phi_{\bm c}[h_1] *_+ \Phi_{\bm c}[h_2])] \in \Lz^1(\pRz^d, \bm x^{\ushort{\bm c-\bm1}})$. The fact that $(h_1*h_2)=\Phi_{\bm c}^{-1}[(\Phi_{\bm c}[h_1] *_+ \Phi_{\bm c}[h_2])] $ follows from simple calculus.
Through this representation we see that 
\begin{align*}
\mathcal M_{\bm c}[(h_1*h_2)]&=(2\pi)^{d/2} \mathcal F[\Phi_{\bm c}[(h_1*h_2)]] = (2\pi)^{d/2} \mathcal F[(\Phi_{\bm c}[h_1]*_+ \Phi_{\bm c}[h_2])] \\
&=(2\pi)^{d}\mathcal F[\Phi_{\bm c}[h_1]] \mathcal F[\Phi_{\bm c}[h_2]]
= \mathcal M_{\bm c}[h_1] \mathcal M_{\bm c}[h_2].
\end{align*}
In analogy, we can see that if additionally $h_1,h_2\in \mathbb L^1(\pRz^d, \bm x^{\ushort{\bm c-\bm 1}})$ and $h_2\in \mathbb L^2(\pRz^d, \bm x^{\ushort{2\bm c- \bm1}})$, then $\Phi_{\bm c}[h_1] *_+ \Phi_{\bm c}[h_2] \in \mathbb L^1(\Rz^d) \cap \mathbb L^2(\Rz^d)$ which implies that $h_1*h_2\in \mathbb L^2(\pRz^d, \bm x^{\ushort{2\bm c-\bm1}})$.

In the upcoming proofs we are in need of the following inequalities. The first inequality is due to
\cite{Talagrand1996}, the formulation can be found for
example in \cite{KleinRio2005}.
\begin{lemma}[Talagrand inequality]\label{tal:re} Let
	$X_1,\dots,X_n$ be independent $\mathcal Z$-valued random variables and let \begin{align*}\bar{\nu}_{h}=n^{-1}\sum_{i=1}^n\left[\nu_{h}(X_i)-\E\left(\nu_{h}(X_i)\right) \right]\end{align*} for $\nu_{h}$ belonging to a countable class $\{\nu_{h},h\in\mathcal H\}$ of measurable functions. Then,
	\begin{align}
	\E(\sup_{h\in\mathcal H}|\overline{\nu}_h|^2-6\Psi^2)_+\leq C \left[\frac{\tau}{n}\exp\left(\frac{-n\Psi^2}{6\tau}\right)+\frac{\psi^2}{n^2}\exp\left(\frac{-K n \Psi}{\psi}\right) \right]\label{tal:re1} 
	\end{align}
	with numerical constants $K=({\sqrt{2}-1})/({21\sqrt{2}})$ and $C>0$ and where
	\begin{equation*}
	\sup_{h\in\mathcal H}\sup_{z\in\mathcal Z}|\nu_{h}(z)|\leq \psi,\qquad \E(\sup_{h\in \mathcal H}|\bar{\nu_{h}}|)\leq \Psi,\qquad \sup_{h\in \mathcal H}\frac{1}{n}\sum_{i=1}^n \Var(\nu_{h}(X_i))\leq \tau.
	\end{equation*}
\end{lemma}
\begin{remark}
	Keeping the bound \eqref{tal:re1}  in mind, let us
	specify particular choices $K$, in fact $K\geq \tfrac{1}{100}$.
	The next bound is now an immediate consequence, 
	\begin{align}
	&\E(\sup_{h\in \mathcal H}|\overline{\nu}_{h}|^2-6\Psi^2)_+\leq C \left(\frac{\tau}{n}\exp\left(\frac{-n\Psi^2}{6\tau}\right)+\frac{\psi^2}{n^2}\exp\left(\frac{-n \Psi}{100\psi}\right) \right).\label{tal:re3} 
	\end{align}
	In the sequel we will make use of the slightly simplified bounds \eqref{tal:re3} rather than \eqref{tal:re1}.
\end{remark}
The next inequality was proven by \cite{Nagaev1979}. A similar formulation can be found in \cite{LiuXiaoWu2013} equation (1.3).

\begin{lemma}[Nagaev's inequality]\label{lem:nagaev} Let $X_1, \dots, X_n$ be i.i.d. mean-zero random variables with $\E(|X_1|^p)< \infty$ for $p>2$. Then for any $x\in \pRz$ holds
	\begin{align*}
	\mathbb P(\sum_{j\in \nset{n}} X_j \geq x)  \leq (1+2p^{-1} )^p \frac{n\E(|X_1|^p)}{x^p} +2\exp(-\frac{a_px^2}{n \E(X_1^2)})
	\end{align*} 
	where $a_p=2e^{-p}(p+2)^{-2}$.
	
\end{lemma}

\subsection{Proofs of Section \ref{mt}}\label{a:mt}
% --------------------------------------------------------------------
% <<Proof of Re key argument>>
% --------------------------------------------------------------------

\begin{proof}[\textbf{\upshape Proof of Proposition \ref{pr:consis}}]
	For $\bm k \in \pRz^d$ we see that $f-f_{\bm k} \in \Lz^2(\pRz^d, \bm x^{\ushort{2\bm c- \bm 1}})$ with $\Mela{f-f_{\bm k}}{\bm c}= \Mela{f}{\bm c} \1_{\pRz^d\setminus Q_{\bm k}}$. We deduce by application of the Plancherel equality that $\langle f-f_{\bm k}, f_{\bm k}- \widehat f_{\bm k} \rangle_{\bm x^{\ushort{2\bm c-\bm 1}}}=\langle \Mela{f}{\bm c} \1_{\pRz^d\setminus Q_{\bm k}}, (\Mela{f}{\bm c}- \widehat{\mathcal M}_{\bm c}) \1_{Q_{\bm k}} \rangle_{\Rz^d} =0$ which implies that
	\begin{align*}
	\E_{f_Y}^n (\|f-\widehat{f}_{\bm k}\|_{\bm x^{\ushort{2\bm c-\bm 1}}}^2) = \|f-f_{\bm k}\|_{\bm x^{\ushort{2\bm c-\bm 1}}}^2 + \E_{f_Y}^n(\|\widehat f_{\bm k}- f_{\bm k} \|_{\bm x^{\ushort{2\bm c-\bm 1}}}^2).
	\end{align*}
	Now by application of the Fubini-Tonelli theorem we interchange the integration order to get
	\begin{align*}
	\E_{f_Y}^n(\|\widehat f_{\bm k}- f_{\bm k} \|_{\bm x^{\ushort{2\bm c-\bm 1}}}^2) &= \frac{1}{(2\pi)^d} \int_{Q_{\bm k}} \E_{f_Y}^n(|\Mela{f_Y}{c}(\bm t)-\widehat{\mathcal M}_{\bm c}(\bm t)|^2) |\Mela{g}{\bm c}(\bm t)|^{-2} d\bm t \leq \frac{1}{(2\pi)^d n} \sigma \Delta_g(\bm k).
	\end{align*}
\end{proof}

\begin{proof}[\textbf{\upshape Proof of Theorem \ref{theorem:lower_bound}}]
	First we outline the main steps of the proof. Let us denote by $\mathcal I:=\{j \in \nset{d}: c_j>1/2\}$ the subset of indices.  We will construct a family of
	functions in $\Dz_{\bm c}^{\bm s}(L)$ by a perturbation of the
	density $f_o: \pRz \rightarrow \pRz$ with small bumps, such that their $\Lz^2(\pRz^d, \bm x^{2\bm c-\bm 1})$-distance
	and the Kullback-Leibler divergence of their induced distributions can be bounded from below and
	above, respectively. The claim then follows then by applying Theorem 2.5
	in \cite{Tsybakov2009}. We use
	the following construction, which we present first.
	
	Denote by $C_c^{\infty}(\Rz)$  the set of all smooth functions with compact support in $\Rz$ and let $\psi\in C_c^{\infty}(\Rz)$ be a function with support in $[0,1]$ and $\int_0^{1} \psi(x)dx = 0$. For each $j\in \nset{d}$ and $K_j\in\Nz$ (to be selected below) and
	$k_j\in\nsetro{0,K_j}$ we define the bump-functions
	$\psi_{k_j, K_j}(x_j):= \psi(x_jK_j-K_j-k_j),$ $x_j\in\Rz$ and define for $p\in \Nz_0:=\{z\in \mathbb Z: z\geq 0\}$ the finite constant $C_{p,\infty}:= \max(\|\psi^{(l)}\|_{\infty}, l\in \nset{0,p})$. Let us further define the operator $\mathrm S: C_c^{\infty}(\Rz)\rightarrow C_c^{\infty}(\Rz)$ with $\mathrm S[f](x)=-x f^{(1)}(x)$ for all $x\in \Rz$ and define $\mathrm S^1:=\mathrm S$ and $\mathrm S^{n}:=\mathrm S \circ \mathrm S^{n-1}$ for $n\in \Nz, n\geq 2$.   Now, for $p \in \Nz$, we define the function $
	\psi_{k_j,K_j,p}(x_j):= \mathrm S^{p} [\psi_{k_j,K_j}](x_j)=\sum_{i=1}^{p} c_{i,p} x_j^i K_j^{i} \psi^{(i)}(x_jK_j-K_j-k_j)$ for $x_j \in \pRz$ and $c_{i,p} \geq 1$ and let $c_p:= \sum_{i=1}^p c_{i,p}$ \\
	For a bump-amplitude $\delta>0, \bm\gamma\in \mathbb N^d$ and $\bm K:=(K_1, \dots, K_d)^T\in \Nz^d$  define 
	\begin{align*}
	\mathcal K:= \bigtimes_{j\in \nset{d}}\nsetro{0,K_j}:=\nsetro{0,K_1} \times \cdots \times \nsetro{0,K_d}=\{\bm k\in \mathbb N_0^d: \forall j\in \nset{d}: k_j < K_j\}.
	\end{align*} and a vector $\bm \theta =( \theta_{\bm k+\bm 1})_{\bm k\in \mathcal K} \in \{0,1\}^{\bigtimes_{j\in \nset{d}} \nset{K_j}}=:\bm \Theta$ we define
	\begin{equation}\label{equation:lobodens}
	\hspace*{-0.5cm}f_{\bm{\theta}}(\bm x)=f_o(\bm x)+ \delta \mathrm F_{\bm K,\bm \gamma, \bm s}^{-1/2} \sum_{\bm k=(k_1,\dots, k_d)^T\in \mathcal K}\theta_{\bm k+\bm 1} \prod_{j\in \nset{d}}
	\psi_{k_j, K_j,\gamma_j}(x_j),
	\end{equation}
	where $\mathrm{F}_{\bm K, \bm \gamma, \bm s}:=\bm K^{\ushort{2\bm \gamma}} \sum_{j\in \nset{d}} K_j^{2s_j}$ and $ f_o(\bm x):=\prod_{j\in \nset{d}} f_{o,j}(x_j)$ with
	\begin{align*}
	f_{o,j}(x):= \begin{cases}
	\exp(-x)\1_{\pRz}(x) ,& j\in \mathcal I; \\
	x\exp(-x)\1_{\pRz}(x), &\text{ else. }
	\end{cases}
	\end{align*}
	Until now, we did not give a sufficient condition to ensure that our constructed functions $\{f_{\bm{\theta}}: \bm{\theta} \in \bm \Theta\}$ are in fact densities. This condition is given by the following lemma.
	\begin{lemma}\label{mm:lem:den}
		Let $0<\delta< \delta_o(\psi,\bm \gamma):=\exp(-2d)/(\prod_{j\in \nset{d}}  2^{\gamma_j}C_{\gamma_j, \infty} c_{\gamma_j})$. Then for all $\bm{\theta}\in\bm \Theta$, $f_{\bm{\theta}}$ is a density.
	\end{lemma}
	Further, one can show that these densities  all lie inside the ellipsoids $\Dz^{\bm s}_{\bm c}(L)$ for $L$ big enough. This is captured in the following lemma.
	\begin{lemma}\label{lemma:Lag_SobDen}Let
		$\bm\wSob\in\Nz^d$. Then,
		there is $L_{\bm s, \bm \gamma,\bm c,\delta}>0$ such that $f_o$
		and any $f_{\bm{\theta}}$ as in  \eqref{equation:lobodens} with
		$\bm\theta\in \bm \Theta$,  belong to $\Dz^{\bm s}_{\bm c}(L_{\bm s, \bm \gamma,\bm c,\delta})$.
	\end{lemma}
	For sake of simplicity we denote for a function $\varphi \in \Lz^2 (\pRz^d, \bm x^{\ushort{2\bm{\tilde c} -\bm 1}})$ the multiplicative convolution with $g$ by $\widetilde{\varphi} := (\varphi *g)$. Futher we see that for  $\bm y\in (0,2)^d$  holds
	\begin{align}\label{eq:dens:low}
	\widetilde{f_o}(\bm y)&=\bm y^{\ushort{\bm 1-2\bm{\tilde c}}}\int_{\pRz^d} g(\bm x)  \bm x^{\ushort{2\bm{\tilde c}}}  \prod_{j\in \nset{d}} \exp(-y_j/x_j)d\bm x \geq \bm y^{\ushort{\bm 1-2\bm{\tilde c}}}\int_{\pRz^d} g(\bm x)  \bm x^{\ushort{2\bm{\tilde c}}}  \prod_{j\in \nset{d}} \exp(-2/x_j)d\bm x=:c_g \bm y^{\ushort{\bm 1-2\bm{\tilde c}}}
	\end{align}  where $c_{g}>0$ since otherwise $g=0$ almost everywhere.
	Exploiting \textit{Varshamov-Gilbert's
		lemma} (see \cite{Tsybakov2009}) in Lemma \ref{lemma:tsyb_vorb} we show
	further that there is $M\in\Nz$ with $M\geq 2^{\prod_{j\in \nset{d}} K_j/8}$ and a subset
	$\{\bm \theta^{(0)}, \dots, \bm \theta^{(M)}\}$ of $\bm \Theta$ with
	$\bm \theta^{(0)}=(0, \dots, 0)$ such that for all
	$j, l \in \nset{0, M}$, $j \neq l$ the $\Lz^2(\pRz^d, \bm x^{2\bm c-\bm 1})$-distance and the
	Kullback-Leibler divergence are bounded for $\bm K\geq \bm K_o(\bm \gamma,\bm c,\psi)$.
	\begin{lemma}\label{lemma:tsyb_vorb}
		Let $\bm K\geq \bm K_o(\psi,\bm\gamma,\bm c)$ understood componentwise). Then there exists a subset $\{\bm \theta^{(0)}, \dots, \bm \theta^{(M)}\}$ of $\bm \Theta$  with $\bm \theta^{(0)}=(0, \dots, 0)$ such that $M\geq 2^{8^{-1}\prod_{j\in \nset{d}} K_j}$ and for all $j, l \in \llbracket 0, M \rrbracket, j \neq l$ holds 
			\begin{enumerate}
		\item[\upshape (i)] $\| f_{\bm\theta^{(j)}}- f_{\bm\theta^{(l)}}\|_{\bm x^{2\bm c-\bm 1}}^2 \geq \frac{C_{\bm \gamma, \bm c}\delta^2}{\sum_{j\in \nset{d}} K_j^{2s_j}},$
		\item[\upshape (ii)]  $\text{KL}(\widetilde{f}_{\bm\theta^{(j)}}, \widetilde{f}_{\bm\theta^{(0)}}) \leq  \frac{ C_{g,\bm \gamma}\delta^2 \log(M) \bm K^{\ushort{-2\bm \gamma-\bm 1}}} {\sum_{j\in \nset{d}} K_j^{2s_j}},$
		\end{enumerate} 
		 where $\text{KL}$ is the Kullback-Leibler-divergence.
	\end{lemma}
	Selecting $K_j=\lceil n^{1/(2s_j+s_j \sum_{i\in \nset{d}} (2\gamma_i+1) s_i^{-1}}\rceil$, it follows that for $n\geq n_{\bm s,\bm \gamma}$
	 \begin{align*}
	(\sum_{j\in \nset{d}} K_j^{2s_j}) ^{-1} \geq c n^{-1/(1+0.5  \sum_{i\in \nset{d}} (2\gamma_i+1) s_i^{-1})}, \quad\frac{\bm K^{\ushort{-2\bm\gamma-\bm 1}}}{\sum_{j\in \nset{d}} K_j^{2s_j}} \leq  n^{-1}
	\end{align*}
	 and thus
	\begin{align*}
	\frac{1}{M}\sum_{j=1}^M
	\text{KL}((\widetilde{f}_{\bm{\theta^{(j)}}})^{\otimes
		n},(\widetilde{f}_{\bm{\theta^{(0)}}})^{\otimes n})
	&= \frac{n}{M} \sum_{j=1}^M \text{KL}(
	\widetilde{f}_{\bm{\theta^{(j)}}},\widetilde{f}_{\bm{\theta^{(0)}}} )
	\leq c_{\delta,g,\bm \gamma}  \log(M),
	\end{align*}
	where $c_{\delta,g,\bm \gamma}< 1/8$ for all
	$\delta\leq \delta_1( g, \bm \gamma,\bm s)$ and $M\geq 2$ for
	$n\geq n_{\bm s,\bm \gamma}$. Thereby, we can use Theorem 2.5 of
	\cite{Tsybakov2009}, which in turn for any estimator $\widehat f$ of $f$
	implies
	\begin{align*}
	\sup_{f\in\Dz^{\bm s}_{\bm c}(L)}
	\mathbb P\big(\|\widehat f-f\|_{\bm x^{2\bm c-\bm 1}}^2\geq
	\frac{c_{\delta, \bm\gamma, \bm c}}{2}n^{-1/(1+0.5  \sum_{i\in \nset{d}} (2\gamma_i+1) s_i^{-1})} \big)\geq 0.07.
	\end{align*}
	Note that the constant $c_{\delta,\bm \gamma, \bm c}$ does only depend on
	$\psi,\gamma $ and $\delta$, hence 
	it is independent of the parameters $\wSob,\rSob$ and $n$. The claim
	of Theorem \ref{theorem:lower_bound} follows by using Markov's inequality,
	which completes the proof.
\end{proof}

\begin{proof}[\textbf{\upshape Proof of Lemma \ref{mm:lem:den}}]
	For any $h\in C_c^{\infty}(\Rz)$ we can state that $\int_{-\infty}^{\infty} \mathrm S[h](x)dx = [-x h(x)]^{\infty}_{-\infty} + \int_{-\infty}^{\infty} h(x) dx=\int_{-\infty}^{\infty} h(x) dx$ and therefore $\int_{-\infty}^{\infty} \mathrm S^p[h](x)dx = \int_{-\infty}^{\infty} h(x)dx$ for $p\in \Nz$. Thus for every $j\in \nset{d}$ we get $\int_{-\infty}^{\infty} \psi_{k_j,K_j,\gamma_j}(x_j)dx_j = \int_{-\infty}^{\infty} \psi_{k_j,K_j}(x_j)dx_j =0$ which implies that for any $\delta >0$ and $\bm{\theta}\in\bm \Theta$ we have $\int_{\pRz^d } f_{\bm{\theta}}(\bm x)d\bm x= 1$.\\
	Now due to the construction \eqref{equation:lobodens} of the functions $\psi_{k_j,K_j}$ we easily see that the function  $\psi_{k_j,K_j}$ has support on $[1+k_j/K_j,1+(k_j+1)/K_j]$ which leads to  $\psi_{k_j,K_j}$ and $ \psi_{l_j,K_j}$ having disjoint supports if $k_j\neq l_j$. Here, we want to emphasize that $\mathrm{supp}(\mathrm S[h]) \subseteq \mathrm{supp}(h)$ for all $h\in C_c^{\infty}(\Rz)$. This implies that $\psi_{k_j,K_j,\gamma_j}$ and $ \psi_{l_j,K_j,\gamma_j}$ have disjoint supports if $k_j\neq l_j$, too.
	For $\bm x\in \pRz^d\setminus \bigtimes_{j\in \nset{d}} [1, 2]$ we have $f_{\bm{\theta}}(\bm x)=\exp(-\sum_{j\in\nset{d}} x_j)\prod_{j \in \mathcal I^c} x_j\geq 0$. Now let us consider the case $\bm x\in\bigtimes_{j\in \nset{d}} [1,2]$. In fact there are $k_{1,o}\in\nsetro{0,K_1}, \dots, k_{d,o} \in \nsetro{0,K_d}$ such that $\bm x \in \bigtimes_{j\in \nset{d}} [1+k_{j,o}/K_j,1+ (k_{j,o}+1)/K_j]$ and hence for $\bm k_o:=(k_{1,o},\dots, k_{d,o})^T$
	\begin{equation*}
	\hspace*{-0.5cm}f_{\bm{\theta}}(\bm x)= f_o(\bm x) + \delta \mathrm F_{\bm K,\bm \gamma, \bm s}^{-1/2} \theta_{\bm k_o+\bm 1}\prod_{j\in \nset{d}}\psi_{k_{o,j},K_j, \gamma_j}(x_j) \geq \exp(-2d)  - \delta \prod_{j\in \nset{d}}  2^{\gamma_j}C_{\gamma_j, \infty} c_{\gamma_j},
	\end{equation*}
	since $\|\psi_{k_j,K_j,\gamma_j}\|_{\infty} \leq 2^{\gamma_j} C_{\gamma_j, \infty} c_{\gamma_j} K_j^{\gamma_j}$ for any $k_j\in \nsetro{0,K_j}$ and $j\in \nset{d}$ and $\mathrm F_{\bm K, \bm \gamma, \bm s}\geq 1$. Choosing $\delta\leq \delta_o(\psi,\bm \gamma)=\exp(-2d)/(\prod_{j\in \nset{d}}  2^{\gamma_j}C_{\gamma_j, \infty} c_{\gamma_j})$ ensures $f_{\bm{\theta}}(x) \geq 0$ for all $x\in \pRz.$ 
\end{proof}

\begin{proof}[\textbf{\upshape Proof of Lemma \ref{lemma:Lag_SobDen}}]Our proof starts with the
	observation that $f_o(\bm x)= \prod_{j\in \nset{d}} f_{o,j}(x_j)$ where $f_{o,j}(x_j):=\exp(-x_j)$ if $j\in \mathcal I$ and $f_{o,j}(x_j):=x_j\exp(-x_j)$ else, for all $\bm x \in \pRz^d$. By the definition of the multivariate Mellin transform, compare \eqref{eq:mel:l1}, we see that $f_{o} \in \Lz^2(\pRz, \bm x^{\ushort{2\bm c-\bm 1}}) \cap \Lz^1(\pRz^d, \bm x^{\ushort{\bm c-\bm 1}})$ holds for every $\bm c \in \pRz^d$  and that
	for all $t\in \Rz^d$ we have 
	\begin{align*}
	\Mela{f_o}{\bm c}(\bm t)= \prod_{j\in \nset{d}} \Mela{f_{o,j}}{c_j}(t_j)= \prod_{j \in \nset{d}} \Gamma(c_j+it_j) \cdot\prod_{j\in \mathcal I^c} (c_j+it_j).
	\end{align*}
	Now by applying the Stirling formula (see also \cite{BelomestnyGoldenshluger2020}) we get $|\Gamma(c_j+it_j)| \sim |t_j|^{c_j+1/2} \exp(-\pi|t_j|/2 )$, $|t|\geq 2$. Thus for every $\bm s\in \Nz^d$ there exists $L_{\bm s, \bm c}$ such that $|f_o|_{\bm s}^2 \leq L $ for all $L\geq L_{\bm s, \bm c}$. 
	
	Next we consider $|f_o-f_{\bm{\theta}}|_{\bm s}$. Again we see that, $f_o- f_{\bm \theta} \in C_c^{\infty}(\pRz^d) \subset \Lz^2(\pRz^d, \bm x^{\ushort{2\bm c- \bm 1}}) \cap \Lz^1(\pRz^d, \bm x^{\ushort{\bm c-\bm1}})$ with 
	\begin{align*}
	\Mela{f_o-f_{\bm \theta}}{\bm c}(\bm t)=\delta \mathrm F_{\bm K, \bm \gamma, \bm s}^{-1/2} \sum_{\bm k\in \mathcal K} \theta_{\bm k+\bm 1}\prod_{j\in \nset{d}} \Mela{\psi_{k_j, K_j,\gamma_j}}{c_j}(t_j), \quad  \bm t\in \Rz^d.
	\end{align*}
	Now for any fixed $\iota\in\nset{d}$ we derive from $\psi_{k_{\iota}, K_{\iota}, \gamma_{\iota} }= \mathrm S^{\gamma_{\iota}} [\psi_{k_{\iota}, K_{\iota}}]$  for any $t_j\in \Rz$ that $\Mela{\psi_{K_{\iota},k_{\iota},\gamma_{\iota}}}{c_j}(t_{\iota})=(c_{\iota}+it_{\iota})^{-s_{\iota}} \Mela{\psi_{K_{\iota},k_{\iota},\gamma_{\iota}+s_{\iota}}}{c_{\iota}}(t_{\iota})$.
	This implies that
	\begin{align*}
	\hspace*{-1cm}|(1+t_{\iota})^{s_{\iota}} &\Mela{f_o-f_{\bm \theta}}{\bm c}(\bm t)|^2\leq C_{c_{\iota}, s_{\iota}} \delta^2 \mathrm F_{\bm K, \bm \gamma, \bm s}^{-1} \left| \sum_{\bm k\in \mathcal K} \theta_{\bm k+\bm 1} \Mela{\psi_{k_{\iota}, K_{\iota}, \gamma_{\iota}+s_{\iota}}}{c_{\iota}}(t_{\iota})\prod_{j\in \nset{d}, j\neq \iota} \Mela{\psi_{k_j, K_j,\gamma_j}}{c_j}(t_j) \right|^2.
	\end{align*}
	Now using that the inverse Mellin operator is linear and by a factorization argument we get 
	\begin{align*}
	\mathcal M^{-1}_{\bm c}[&\sum_{\bm k\in \mathcal K} \theta_{\bm k+\bm 1} \Mela{\psi_{k_{\iota}, K_{\iota}, \gamma_{\iota}+s_{\iota}}}{c_{\iota}}(t_{\iota})\prod_{j\in \nset{d}, j\neq \iota} \Mela{\psi_{k_j, K_j,\gamma_j}}{c_j}(t_j)](\bm x) = \sum_{\bm k\in \mathcal K} \theta_{\bm k+\bm 1} \psi_{k_{\iota}, K_{\iota}, \gamma_{\iota}+s_{\iota}}(x_{\iota}) \prod_{j\in \nset{d}, j\neq \iota} \psi_{k_j, K_j,\gamma_j}(x_j)
	\end{align*}
	which implies that due to the disjoint supports of $\psi_{k_j, K_j, \gamma_j}$ and another factorization argument,
	\begin{align*}
	\hspace*{-1cm}\|(1+t_{\iota})^{s_{\iota}} \Mela{f_o-f_{\bm \theta}}{\bm c}\|^2_{\Rz^d} &\leq C_{\bm c_{\iota}, s_{\iota}} \delta^2 \mathrm F_{\bm K, \bm \gamma, \bm s}^{-1}\sum_{\bm k\in \mathcal K} \|\psi_{k_{\iota}, K_{\iota}, \gamma_{\iota}+s_{\iota}}\|_{ x^{2 c_j- 1}}^2 \prod_{j\in \nset{d}, j\neq \iota} \| \psi_{k_j, K_j,\gamma_j}\|_{ x^{2 c_j- 1}}^2 \leq C_{\bm c, \bm s, \bm \gamma} \delta^2 \frac{K_{\iota}^{2s_{\iota}}}{\sum_{j\in \nset{d}} K_j^{2s_j}}
	\end{align*}
	since $\| \psi_{k_j, K_j,p}\|_{ x^{2 c_j- 1}}^2 =\int_{\pRz} \psi_{k_j, K_j,p}^2(x) x^{2c_j-1}dx \leq C_{\gamma_j, s_j,c_j} K_j^{2p-1}$ for any $p\in \Nz$.
	We follow $|f_o-f_{\bm \theta}|_{\bm s,\bm c}^2 \leq C_{\bm c, \bm s, \bm \gamma}$.  Finally, we have that  $|f_{\bm{\theta}}|_{\bm s}^2 \leq 2(|f_o-f_{\bm{\theta}}|_{\bm s}^2 + |f_o|_{\bm s}^2) \leq 2(C_{(\bm s, \bm \gamma,\bm c,\delta,\psi)}+ L_{\bm s}) =: L_{\bm s, \bm \gamma,\bm c,\delta,1}$.
	
	Now we have to consider the moment condition $\E_{f_{\bm \theta}}(\bm X^{\ushort{2\bm  c-\bm 2}}) \leq L$. In fact we have
	\begin{align*}
	\int_{\pRz^d} \bm x^{\ushort{2\bm c-\bm 2}} f_{\bm \theta}(\bm x) d\bm x &= \prod_{j\in \nset{d}} \int_{\pRz} x^{2c_j-2} f_{o,j}(x_j) dx_j + \delta \mathrm F_{\bm K, \bm\gamma, \bm s}^{-1/2} \sum_{\bm k=(k_1,\dots, k_d)^T\in \mathcal K} \theta_{\bm k+\bm 1} \prod_{j\in \nset{d}} \int_{\pRz} x^{2c_j-2} \psi_{k_j, K_j, \gamma_j}(x_j) dx_j\\
	&\leq C_{\bm c} +\delta  \mathrm F_{\bm K, \bm\gamma, \bm s}^{-1/2} \sum_{\bm k=(k_1,\dots, k_d)^T\in \mathcal K} \prod_{j\in \nset{d}} C_{\bm \gamma, \bm c} K_j^{\gamma-1} \leq C_{\bm c} + \delta C_{\bm \gamma, \bm c} =: L_{\bm \gamma, \bm c, \delta, 2}.
	\end{align*}
	Now we choose $L_{\bm s, \bm \gamma, \bm c, \delta}:= \max(L_{\bm s, \bm \gamma, \bm c, \delta, 1}, L_{\bm \gamma, \bm c, \delta, 2})$.
\end{proof}
\begin{proof}[\textbf{\upshape Proof of Lemma \ref{lemma:tsyb_vorb}}]  (i):
	Using  that the functions $(\psi_{k_j,K_j, \gamma_j})$ with different index $k_j$ have disjoint supports and a factorization argument we get 
	\begin{align*}
	\| f_{\bm{\theta}}-f_{\bm{\theta}'}\|_{\bm x^{2\bm c-\bm 1}}^2&= \delta^2 \mathrm F_{\bm K, \bm \gamma, \bm s}^{-1} \| \sum_{\bm k \in \mathcal K} ( \theta_{\bm k+\bm 1}-\theta'_{\bm k+\bm 1}) \prod_{j\in \nset{d}} \psi_{k_j,K_j, \gamma_j}(x_j)\|_{\bm x^{\ushort{2\bm c-\bm 1}}}^2 = \delta^2   F_{\bm K, \bm \gamma, \bm s}^{-1} \sum_{\bm k\in \mathcal K}  ( \theta_{\bm k+\bm 1}-\theta'_{\bm k+\bm 1})^2 \prod_{j\in \nset{d}} \|   \psi_{k_j,K_j, \gamma_j}\|_{x^{2c_j-1}}^2 \\
	&\geq \delta^2   F_{\bm K, \bm \gamma, \bm s}^{-1} \rho( \theta_{\bm k+\bm 1},\theta'_{\bm k+\bm 1})^2  C_{\bm \gamma, \bm c} \bm K^{\ushort{2\bm \gamma -\bm 1}},
	\end{align*}
	where the last step follows if we can show that there exists a $c_{c_j}>0$ 
	\begin{align}\label{eq:one:dim}
	\|\psi_{k_j,K_j, \gamma_j}\|_{x^{2c_j-1}}^2 = \int_0^{\infty} \psi_{k_j, K_j,\gamma_j}(x)^2 x^{2c_1-1} dx \geq \frac{c_{c_j}K_j^{2\gamma_j-1}\|\psi^{(\gamma_j)}\|_{x^0}^2 }{2} \
	\end{align}  for $K_j$ big enough. Here $\rho(\bm \theta, \bm \theta'):= \sum_{\bm k\in \mathcal K} \1_{\{\bm \theta_{\bm k+1} \neq \bm \theta'_{\bm k+1}\}}$ denotes the Hamming distance.\\
	To show \eqref{eq:one:dim} we observe that 
	\begin{align*}
	\| \psi_{k_j,K_j, \gamma_j} \|_{x^{2c_j+1}}^2 = \sum_{i,\iota \in \llbracket 1, \gamma_j\rrbracket} c_{i,\gamma_j}c_{\iota,\gamma_j} \int_0^{\infty} x^{\iota+i+1} \psi_{k_j,K_j}^{(\iota)}(x)  \psi_{k_j,K_j}^{(i)}(x)dx
	\end{align*}   and by defining $\Sigma:= \| \psi_{k_j,K_j,\gamma_j}\|_{x^{2c_j+1}}^2 -  \int_0^{\infty}(x^{\gamma_j} \psi_{k_j,K_j}^{(\gamma_j)}(x))^2 xdx$ we can show
	\begin{align}\label{equation:low_l2_gamma}
	\| \psi_{k_j,K_j, \gamma_j} \|_{x^{2c_j-1}}^2 = \Sigma + \int_0^{\infty}(x^{\gamma} \psi_{0,K}^{(\gamma)}(x))^2 x^{2c_j-1}dx \geq \Sigma + c_{c_j}K^{2\gamma_j-1} \| \psi^{(\gamma_j)} \|^2 \geq \frac{c_{c_j} K_j^{2\gamma_j-1} \| \psi^{(\gamma_j)} \|^2}{2}
	\end{align}
	as soon as $|\Sigma|\leq \frac{ c_{c_j}K^{2\gamma-1} \| \psi^{(\gamma_j)} \|^2}{2}$. This is obviously true as soon as $K_j \geq K_{o}(\gamma_j, c_j, \psi)$ and thus $	\| f_{\bm{\theta}}-f_{\bm{\theta}'}\|_{\bm x^{2\bm c-\bm 1}}^2 \geq  \delta^2C_{\bm \gamma, \bm c}\bm K^{\ushort{-\bm 1}} (\sum_{j\in\nset{d}} K_j^{2s_j} )^{-1}\rho(\bm{\theta},\bm{\theta}')$ for $\bm K\geq \bm K_o(\psi, \bm\gamma,\bm c)$ which is understood in a componentwise sense.
	
	Now let us interpretate the objects $\bm \theta\in \bm \Theta$ as vectors using the canonical bijection $\mathrm T:\bm \Theta\rightarrow\{0,1\}^{\prod_{j\in \nset{d}} K_j}$. Then we have  $\rho(\bm \theta_{\bm k+\bm 1}, \bm \theta'_{\bm k+\bm 1})) = \tilde{\rho}(T(\bm \theta_{\bm k+\bm 1}), T(\bm \theta'_{\bm k+\bm 1}))$ where $\tilde{\rho}(\bm \vartheta, \bm \vartheta') :=\sum_{k=1}^{\prod_{j\in \nset{d}} K_j}  \1_{\{\vartheta_{j}\neq\vartheta'_{j}\}}$ for any $\bm \vartheta, \bm \vartheta'\in \{0,1\}^{\prod_{j\in \nset{d}} K_j}$.
	Using the Varshamov-Gilbert Lemma (see \cite{Tsybakov2009}) which states that for $\prod_{j\in \nset{d}} K_j \geq 8$ there exists a subset $\{\bm \vartheta^{(0)}, \dots, \bm\vartheta^{(M)}\}$ of $\{0,1\}^{\prod_{j\in \nset{d}} K_j}$ with $\bm \vartheta^{(0)}=(0, \dots, 0)$ such that $\tilde \rho(\bm \vartheta^{(j)}, \bm \vartheta^{(k)}) \geq \prod_{j\in \nset{d}} K_j/8$ for all $j ,k\in \llbracket 0, M \rrbracket, j\neq k $ and $M\geq 2^{8^{-1}\prod_{j\in \nset{d}} K_j}$. Defining $\bm \theta^{(j)}:= T^{-1}(\bm\vartheta^{(j)})$ for $j\in \nset{0,M}$  leads to $\| f_{\bm\theta^{(j)}}- f_{\bm\theta^{(l)}}\|_{\bm x^{2\bm c-\bm 1}}^2 \geq C_{\bm \gamma, \bm c}\delta^2(\sum_{j\in \nset{d}} K_j^{2s_j})^{-1}$.
	
	(ii): For the second part we have $f_o=f_{\bm\theta^{(0)}}$, and by using $\text{KL}(\widetilde{f}_{\bm \theta}, \widetilde{f}_o) \leq \chi^2(\widetilde{f}_{\bm\theta},\widetilde{f}_o):= \int_{\pRz^d} (\widetilde{f}_{\bm\theta}(x) -\widetilde{f}_o(\bm x))^2/\widetilde{f}_o(\bm x) d\bm x$ it is sufficient to bound the $\chi$-squared divergence. We notice that since $U_1,\dots, U_n$ are independent we can write $g(\bm x)= \prod_{j\in \nset{d}} g_j(x_j)$ for $\bm x\in\pRz^d$. Further, $\widetilde{f}_{\bm\theta} -\widetilde{f}_o$ has support in $[0,2]^d$ since $f_{\bm{\theta}}-f_o$ has support in $[1,2]^d$ and $g$ has support in $[0,1]^d$. In fact for $\bm y\in \pRz^d$ with $ y_j>2$ for $j\in\nset{d}$,  
	\begin{align*}
	\widetilde{f}_{\bm\theta}(\bm  y) -\widetilde{f}_o(\bm y)&=\int_{\pRz^d} (f_{\bm{\theta}}-f_o)(\bm x)\bm x^{\ushort{-\bm 1}} g(\bm y/\bm x)d\bm x= \delta \mathrm F_{\bm K,\bm \gamma, \bm s}^{-1/2}\sum_{\bm k\in \mathcal K} \bm \theta_{\bm k+\bm 1} \prod_{j\in \nset{d}}  \int_{y_j}^{\infty} \psi_{k_j,K_j,\gamma_j}(x_j) g_j(y_j/x_j) x_j^{-1} d x_j=0.
	\end{align*} 
	Next we have for any $\bm t\in \Rz^d$ by application of assumption of Theorem \ref{theorem:lower_bound}, the convolution theorem and the fact that $\Mela{\psi_{k_j,K_j,\gamma_j}}{\tilde c_j}(t_j)=(\tilde c_j+it_j)^{\gamma_j}\Mela{\psi_{k_j,K_j,0}}{\tilde c_j}(t_j)$
	\begin{align*}
	|\Mela{	\widetilde{f}_{\bm\theta}-\widetilde{f}_o}{\bm{\tilde c}}(\bm t)|&=|\delta \mathrm F_{\bm K,\bm \gamma, \bm s}^{-1/2}\sum_{\bm k\in \mathcal K} \bm \theta_{\bm k+\bm 1}\prod_{j\in \nset {d}} \Mela{\psi_{k_j,K_j,\gamma_j}}{\tilde c_j}(t_j) \Mela{g_j}{\tilde c_j}(t_j)| \leq C_{g,\bm \gamma}\delta \mathrm F_{\bm K,\bm \gamma, \bm s}^{-1/2}|\sum_{\bm k\in \mathcal K} \bm \theta_{\bm k+\bm 1}\prod_{j\in \nset {d}} \Mela{\psi_{k_j,K_j,0}}{\tilde c_j}(t_j)|,
	\end{align*}
	for all $\bm t \in \Rz^d$. 
	Applying the Parseval equality and using the disjoints supports, the factorization property and eq. \ref{eq:dens:low} we get
	\begin{align*}
	\chi^2(\widetilde{f}_{\bm\theta},\widetilde{f}_o)
	\leq C_g \|\widetilde{f}_{\bm \theta}- \widetilde{f}_o\|_{\bm x^{\ushort{2\bm{\tilde c}- \bm 1}}}^2 \leq C_{g,\bm \gamma}  \delta^2 \mathrm F_{\bm K, \bm \gamma, \bm s}^{-1}  \|\sum_{\bm k\in \mathcal K} \prod_{j\in \nset{d}}  \Mela{\psi_{k_j,K_j,0}}{\tilde c_j}(t_j)\|_{\Rz^d}^2 \leq C_{g,\bm \gamma}\delta^2 \mathrm F_{\bm K, \bm \gamma, \bm s}^{-1}   \sum_{\bm k\in \mathcal K} \prod_{j\in \nset{d}} \| \psi_{k_j,K_j,0}\|_{ x^{2\tilde c_j-1}}^2.
	\end{align*}
	The inequality  $\| \psi_{k_j,K_j,0}\|_{ x^{ 2\tilde c_j-1}}^2$. In fact, using that $\|\psi_{k_j, K_j, 0}\|_{x^{2\tilde c_j-1}}^2 \leq K_j^{-1} \|\psi\|_{x^0}^2$ implies $
	\chi^2(\widetilde{f}_{\bm\theta},\widetilde{f}_o) \leq C_{g} \delta^2 \mathrm F_{\bm K, \bm \gamma, \bm s}^{-1} $. 
	Since $M \geq 2^{\prod_{j\in \nset{d}} K_j}$ we can deduce that
	$\mathrm{KL}(\widetilde{f}_{\bm\theta^{(j)}},\widetilde{f}_{\bm\theta^{(0)}})\leq  C_{g,\bm \gamma}\delta^2 \log(M) \bm K^{\ushort{2\bm \gamma-\bm 1}} (\sum_{j\in \nset{d}} K_j^{2s_j})^{-1}.$ 
\end{proof}

\label{a:dd}

\begin{proof}[\textbf{\upshape Proof of Theorem\ref{theo:adap:aniso}}]
	Let $\bm k\in \mathcal K_n$. By definition of the estimator, \eqref{eq:est:def}, we have $Q_{\bm k'}=\mathrm{supp}(\mathcal M_{\bm c}[f_{\bm k'}])$, for $\bm k'\in \mathcal K_n$ and we can find a $\bm K_n\in (\mathbb N^*)^2$ such that for all $\bm k'\in \mathcal K_n$ holds $Q_{\bm k'}\subseteq Q_{\bm K_n}$. Then we have for any $\bm k'\in \mathcal K_n$ that $\|\widehat f_{\bm K_n}\|_{\bm x^{\ushort{2\bm c-\bm 1}}}^2- \|\widehat f_{\bm k'}\|_{\bm x^{\ushort{2\bm c-\bm 1}}}^2=\|\widehat f_{\bm K_n}-\widehat f_{\bm k'}\|_{\bm x^{\ushort{2\bm c-\bm 1}}}^2$ implying with \eqref{eq:data:driven}
	\begin{align}\label{eq:adap:new} \|\widehat f_{\widehat{\bm k}}-\widehat f_{\bm K_n}\|_{\bm x^{\ushort{2\bm c-\bm 1}}}^2 + \widehat{\mathrm{pen}}(\widehat {\bm k}) \leq \|\widehat f_{\bm k}-\widehat f_{\bm K_n}\|^2_{\bm x^{\ushort{2\bm c-\bm 1}}}+ \widehat{\mathrm{pen}}(\bm k).\end{align}
	Now for every $\bm k'\in \mathcal K_n$ we have 
	$ \|\widehat f_{\bm k'}-f_{\bm K_n}\|_{\bm x^{\ushort{2\bm c-\bm 1}}}^2 = \|\widehat f_{\bm k'}-\widehat f_{\bm K_n}\|_{\bm x^{\ushort{2\bm c-\bm 1}}}^2 + \|\widehat f_{\bm K_n}-f_{\bm K_n}\|_{\bm x^{\ushort{2\bm c-\bm 1}}}^2+2\langle \widehat f_{\bm k'}-\widehat f_{\bm K_n}, \widehat f_{\bm K_n}-f_{\bm K_n} \rangle_{\bm x^{\ushort{2\bm c-\bm 1}}}$
	which combined with \eqref{eq:adap:new} implies
	\begin{align} \label{eq:decomp:1:2}
	\| \widehat f_{\widehat{\bm k}}-f_{\bm K_n}\|_{\bm x^{\ushort{2\bm c-\bm 1}}}^2-\|\widehat f_{\bm k}-f_{\bm K_n}\|_{\bm x^{\ushort{2\bm c-\bm 1}}}^2 
	&\leq \widehat{\mathrm{pen}}(\bm k)- \widehat{\mathrm{pen}}(\widehat{\bm k}) + 2\langle \widehat f_{\widehat {\bm k}}-\widehat f_{\bm k}, \widehat f_{\bm K_n}-f_{\bm K_n} \rangle_{\bm x^{\ushort{2\bm c-\bm 1}}}. 
	\end{align}
	Since $
	\langle \widehat f_{\widehat{\bm k}}-\widehat f_{\bm k}, \widehat f_{\bm K_n}-f_{\bm K_n}\rangle_{\bm x^{\ushort{2\bm c-\bm 1}}}= \|\widehat f_{\widehat{\bm k}}-f_{\widehat{\bm k}}\|_{\bm x^{\ushort{2\bm c-\bm 1}}}^2+ \langle f_{\widehat{\bm k}}-f_{\bm k}, \widehat f_{\bm K_n}-f_{\bm K_n}\rangle_{\bm x^{\ushort{2\bm c-\bm 1}}} - \|\widehat f_{\bm k}-f_{\bm k}\|_{\bm x^{\ushort{2\bm c-\bm 1}}}^2
	$
	we get
	\begin{align}\label{eq:decomp:3:2}
	\| \widehat f_{\widehat{\bm k}}-f_{\bm K_n}\|_{\bm x^{\ushort{2\bm c-\bm 1}}}^2 &\leq \| f_{\bm k}-f_{\bm K_n}\|_{\bm x^{\ushort{2\bm c-\bm 1}}}^2 - \|\widehat f_{\bm k}-f_{\bm k}\|_{\bm x^{\ushort{2\bm c-\bm 1}}}^2+2\langle f_{\widehat{\bm k}}-f_{\bm k}, \widehat f_{\bm K_n}-f_{\bm K_n}\rangle_{\bm x^{\ushort{2\bm c-\bm 1}}}+ \widehat{\mathrm{pen}}(\bm k)+ 2\|\widehat f_{\widehat {\bm k}}-f_{\widehat{\bm k}}\|_{\bm x^{\ushort{2\bm c-\bm 1}}}^2-\widehat{\mathrm{pen}}(\widehat{\bm k}).
	\end{align}
	We now consider the term $|2\langle f_{\widehat{\bm k}}-f_{\bm k}, \widehat f_{\bm K_n}-f_{\bm K_n}\rangle_{\bm x^{\ushort{2\bm c-\bm 1}}}|$. First we remind that for any $\bm k'\in \mathcal K_n$
	\begin{align*}
	\|\widehat f_{\bm k'}-f_{\bm k'}\|_{\bm x^{\ushort{2\bm c-\bm 1}}}^2 = \frac{1}{(2\pi)^d}\int_{\mathbb R^d} \mathds 1_{Q_{\bm k'}}(\bm t)\frac{ |\mathcal M_{\bm c}[f_{\bm Y}](\bm t)-\widehat{\mathcal M}_{\bm c}(\bm t)|^2}{|\mathcal M_{\bm c}[g](\bm t)|^2}d\bm t.
	\end{align*}
	Setting  $Q^*:=Q_{\widehat{\bm k}}\cup Q_{\bm k}$ we have $\mathcal M_{\bm c}[f_{\widehat{\bm k}}-f_{ \bm k}]=\mathcal M_{\bm c}[f](\mathds 1_{Q_{\widehat{\bm k}}}- \mathds 1_{Q_{\bm k}})$ implying that $\mathrm{supp}(\mathcal M_{\bm c}[f_{\widehat {\bm k}}-f_{\bm k}]) \subseteq Q^*\subseteq Q_{\bm K_n}$ by definition of $\bm K_n$. Using that $2ab\leq a^2+b^2$ we deduce
	\begin{align*}
	|2\langle f_{\widehat{\bm k}}-f_{\bm k}, \widehat f_{\bm K_n}-f_{\bm K_n}\rangle_{\bm x^{\ushort{2\bm c-\bm 1}}}| &= \frac{2}{(2\pi)^d} \left|\int_{Q^*} \mathcal M_{\bm c}[f_{\widehat{\bm k}}-f_{\bm k}](\bm t) \frac{\widehat{\mathcal M}_{\bm c}(-\bm t)-\mathcal M_{\bm c}[f_{\bm Y}](-\bm  t)}{\mathcal M_{\bm c}[g](-\bm t)} d\bm t \right| \\
	&\leq  \frac{1}{4}\|f_{\widehat{\bm k}}-f_{\bm k}\|_{\bm x^{\ushort{2\bm c-\bm 1}}}^2 + \frac{4}{(2\pi)^d} \int_{\mathbb R^d} \mathds 1_{Q^*}(\bm t)\frac{|\widehat {\mathcal M}_{\bm c}(\bm t)-\mathcal M_{\bm c}[f_{\bm Y}](\bm t)|^2}{|\mathcal M_{\bm c}[g](\bm t)|^2} d\bm t \\
	&\leq \frac{1}{2}\|f_{\widehat{\bm k}}-f_{\bm K_n}\|_{\bm x^{\ushort{2\bm c-\bm 1}}}^2+ \frac{1}{2}\|f_{\bm k}-f_{\bm K_n}\|_{\bm x^{\ushort{2\bm c-\bm 1}}}^2+4\|\widehat f_{\bm k}-f_{\bm k}\|_{\bm x^{\ushort{2\bm c-\bm 1}}}^2 + 4\|\widehat f_{\widehat {\bm k}}-f_{\widehat{\bm k}}\|_{\bm x^{\ushort{2\bm c-\bm 1}}}^2
	\intertext{ using that $\mathds 1_{Q^*} \leq \mathds 1_{Q_k}+ \mathds 1_{Q_{\widehat k}}$. Thus we get }
	|2\langle f_{\widehat{\bm k}}-f_{\bm k}, \widehat f_{\bm K_n}-f_{\bm K_n}\rangle_{\bm x^{\ushort{2\bm c-\bm 1}}}| &\leq   \frac{1}{2} \|\widehat f_{\widehat {\bm k}}-f_{\bm K_n}\|_{\bm x^{\ushort{2\bm c-\bm 1}}}^2+ \frac{1}{2} \|f_{\bm K_n}- f_{\bm k}\|_{\bm x^{\ushort{2\bm c-\bm 1}}}^2+4\|\widehat f_{\bm k}-f_{\bm k}\|_{\bm x^{\ushort{2\bm c-\bm 1}}}^2+\frac{7}{2}  \|\widehat f_{\widehat{\bm k}}-f_{\widehat{\bm k}}\|_{\bm x^{\ushort{2\bm c-\bm 1}}}^2 
	\end{align*}
	implying that
		\begin{align}\label{eq:decomp:3}
	\| \widehat f_{\widehat{\bm k}}-f_{\bm K_n}\|_{\bm x^{\ushort{2\bm c-\bm 1}}}^2 &\leq 3\| f_{\bm k}-f_{\bm K_n}\|_{\bm x^{\ushort{2\bm c-\bm 1}}}^2 +6 \|\widehat f_{\bm k}-f_{\bm k}\|_{\bm x^{\ushort{2\bm c-\bm 1}}}^2+ 2\widehat{\mathrm{pen}}(\bm k)+ 11\|\widehat f_{\widehat {\bm k}}-f_{\widehat{\bm k}}\|_{\bm x^{\ushort{2\bm c-\bm 1}}}^2-2\widehat{\mathrm{pen}}(\widehat{\bm k}).
	\end{align}
	Now since $\E_{f_{\bm Y}}^n(\widehat{\mathrm{pen}}(\bm k))=\mathrm{pen}(\bm k)$ and $6\E_{f_{\bm Y}}^n(\|\widehat f_{\bm k}-f_{\bm k}\|_{\bm x^{\ushort{2\bm c-\bm 1}}}^2) \leq 6\chi^{-1} \mathrm{pen}(\bm k) \leq \mathrm{pen}(\bm k)$ we get combined with \eqref{eq:decomp:3}
	\begin{align*}
	\E_{f_{\bm Y}}^n(\| \widehat f_{\widehat{\bm k}}&-f_{\bm K_n}\|_{\bm x^{\ushort{2\bm c-\bm 1}}}^2) \leq 3\left(\|f_{\bm K_n}-f_{\bm k}\|_{\bm x^{\ushort{2\bm c-\bm 1}}}^2 +\mathrm{pen}(\bm k)\right)+ 11\E_{f_{\bm Y}}^n(\left(\|\widehat f_{\widehat{\bm k}}-f_{\widehat{\bm k}}\|_{\bm x^{\ushort{2\bm c-\bm 1}}}^2-\frac{1}{6}\widehat{\mathrm{pen}}(\widehat{\bm  k})\right)_+) \\
	&\leq 3\left(\|f_{\bm K_n}-f_{\bm k}\|_{\bm x^{\ushort{2\bm c-\bm 1}}}^2 +\mathrm{pen}(\bm k)\right)  + 11\E_{f_{\bm Y}}^n(\left(\|\widehat f_{\widehat{\bm k}}-f_{\widehat{\bm k}}\|_{\bm x^{\ushort{2\bm c-\bm 1}}}^2-\frac{1}{12}\mathrm{pen}(\widehat{\bm  k})\right)_+) + \E_{f_{\bm Y}}^n((\mathrm{pen}(\widehat{\bm k})- 2\widehat{\mathrm{pen}}(\widehat {\bm k}))_+)
	\end{align*}
	The two expectations on the right hand side of the last inequality can be bounded using the following Lemma.
\begin{lemma}\label{lem:talagrand}
	Assume that $\E_{f_Y}(\bm Y^{\ushort{7(\bm c-\bm 1)}})<\infty$ and $\| f_Y \bm x^{\ushort{2\bm c-\bm 1}}\|_{\infty}< \infty$. Then 
	\begin{enumerate}
		\item[\upshape (i)] $\E_{f_{\bm Y}}^n(\left(\|\widehat f_{\widehat{\bm k}} -f_{\widehat{\bm k}} \|^2_{\bm x^{\ushort{2\bm c-\bm 1}}} - \frac{1}{12}\mathrm{pen}(\widehat{\bm k}) \right)_+) \leq C_{f,g,\sigma, \E_{f_Y}(\bm Y^{\ushort{5(\bm c-\bm 1)}})}n^{-1},$
		\item[\upshape (ii)] $\E_{f_{\bm Y}}^n((\mathrm{pen}(\widehat{\bm k})- 2\widehat{\mathrm{pen}}(\widehat {\bm k}))_+)\leq  C_{\sigma, \E_{f_Y}(\bm Y^{\ushort{4(\bm c-\bm 1)}})} n^{-1}.$
	\end{enumerate}
\end{lemma}
	Consequently, we have
	\begin{align*}
	\E_{f_{\bm Y}}^n(\|\widehat f_{\widehat{\bm k}}-f\|_{\bm x^{\ushort{2\bm c-\bm 1}}})^2 ) \leq \|f-f_{\bm K_n}\|_{\bm x^{\ushort{2\bm c-\bm 1}}}^2  + 3\left(\|f_{\bm K_n}-f_{\bm k}\|_{\bm x^{\ushort{2\bm c-\bm 1}}}^2 +\mathrm{pen}(\bm k)\right)  + \frac{C_{f,g}}{n} \leq 3\left(\|f-f_{\bm k}\|_{\bm x^{\ushort{2\bm c-\bm 1}}}^2 +\mathrm{pen}(\bm k)\right) + \frac{C_{f,g}}{n}.
	\end{align*}
	Taking now the infimum over all $\bm k\in \mathcal K_n$ implies the claim.
\end{proof}

\begin{proof}[\textbf{\upshape Proof of Lemma \ref{lem:talagrand}}]
	We start by proving (i). Let us therefore define the set $\Uz:=\{h\in \Lz^2(\pRz^d, \bm x^{\ushort{2\bm c-\bm 1}}): \|h\|_{\bm x^{\ushort{2\bm c-\bm 1}}} \leq 1\}$. Then  for $\bm k'\in \pRz^d$, $\|\widehat f_{\bm k'}-f_{\bm k'}\|_{\bm x^{\ushort{2\bm c-\bm 1}}} = \sup_{h\in \Uz} \langle \widehat f_{\bm k'}-f_{\bm k'}, h \rangle_{\bm x^{\ushort{2\bm c-\bm1}}}$ where 
	\begin{align*}
	\langle \widehat f_{\bm k'}-f_{\bm k'}, h\rangle_{\bm x^{\ushort{2\bm c-\bm 1}}}=  (2\pi)^{-d}  \int_{Q_{\bm k'}} (\widehat{\mathcal M}_{\bm c}(\bm t)- \E_{f_Y}^n(\widehat{\mathcal M}_{\bm c}(\bm t)) \frac{\Mela{h}{\bm c}(-\bm t)}{\Mela{g}{\bm c}(\bm t)} d\bm t,
	\end{align*}
	by application of the Plancherel equality. Now for a sequence $(c_n)_{n\in \Nz}$ we decompose the estimator $\widehat{\mathcal M}_{\bm c}(\bm t)$ into
	\begin{align*}
	\widehat{\mathcal M}_{\bm c}(\bm t):&= n^{-1} \sum_{j\in \nset{n}} \bm Y_j^{\ushort{\bm c-\bm 1+i\bm t}} \1_{(0, c_n)}(\bm Y_j^{\ushort{\bm c-\bm 1}})+ n^{-1} \sum_{j\in \nset{n}} \bm Y_j^{\ushort{\bm c-\bm 1+i\bm t}} \1_{[c_n, \infty)}(\bm Y_j^{\ushort{\bm c-\bm 1}})=: \widehat{\mathcal M}_{\bm c, 1}(\bm t)+\widehat{\mathcal M}_{\bm c, 2}(\bm t)
	\end{align*}
	where $\bm{(0, c_n)}:= (0, c_n)^d$. Setting 
	\begin{align*}
	\nu_{\bm k', i}(h):= \frac{1}{(2\pi)^d} \int_{Q_{\bm k'}} (\widehat{\mathcal M}_{\bm c, i}(\bm t)- \E_{f_Y}^n(\widehat{\mathcal M}_{\bm c, i}(\bm t)) \frac{\Mela{h}{\bm c}(-\bm t)}{\Mela{g}{\bm c}(\bm t)} d\bm t, \quad h\in \Uz, i\in \{1,2\},
	\end{align*} we can deduce that 
	\begin{align}\label{eq:two:process}
\hspace*{-0.7cm}\E_{f_Y}^n(\left(\|\widehat f_{\widehat{\bm k}} -f_{\widehat{\bm k}} \|^2_{\bm x^{\ushort{2\bm c-\bm 1}}} - \frac{1}{12}\mathrm{pen}(\widehat{\bm k}) \right)_+)  &\leq 2\E_{f_Y}^n(\left( \sup_{h\in \Uz} \nu_{\widehat{\bm k}, 1}(h)^2- \frac{1}{24} \mathrm{pen}(\widehat{\bm k}) \right)_+) + 2 \E_{f_Y}^n(\sup_{h\in \Uz}\nu_{\widehat{\bm k},2}(h)^2)).
	\end{align}
	We start by bounding the first summand. To do so, we see that 
	\begin{align*}
	\E_{f_Y}^n( \left( \sup_{h\in \Uz} \nu_{\widehat{\bm k}, 1}(h)^2- \frac{1}{24} \mathrm{pen}(\widehat{\bm k}) \right)_+)  \leq \sum_{\bm k'\in \mathcal K_n}\E_{f_Y}^n(\left( \sup_{h\in \Uz} \nu_{\bm k', 1}(h)^2- \frac{1}{24} \mathrm{pen}(\bm k') \right)_+). 
	\end{align*}
	To control each summand we apply the Talagrand inequality, see Remark \ref{tal:re}, which can be done since there exists a dense subset of $\Uz$. For each $\bm k'\in \pRz^d, h\in \Uz$ we set
	\begin{align*}
	\nu_h(\bm y):=\frac{1}{(2\pi)^d} \int_{Q_{\bm k'}} \bm y^{\ushort{\bm c-\bm 1+i\bm t}} \1_{(0,c_n)}(\bm y^{\ushort{\bm c-\bm 1}}) \frac{\Mela{h}{\bm c}(-\bm t)}{\Mela{g}{\bm c}(\bm t)} d\bm t, \quad \bm y\in \pRz^d.
	\end{align*}
	So,  $\overline\nu_h= \nu_{\bm k',1}(h)$ in the notation of Remark \ref{tal:re}. Thus we need to determine the parameters $\tau, \Psi^2, \psi$. Let us begin with $\Psi$. For $h\in \Uz$ we have  $1\geq\|h\|_{\bm x^{2\bm c-1}}^2=(2\pi)^{-d}\|\Mela{h}{c}\|_{\Rz^d}$. Using the Cauchy Schwartz inequality delivers
	\begin{align*}
	\E_{f_Y}(\sup_{h\in \Uz}\overline{\nu}_h^2)
	&\leq  (2\pi)^{-d} \int_{Q_{\bm k'}} \E_{f_Y}^n(|\widehat{\mathcal M}_{\bm c, 1}(\bm t)- \E_{f_Y}^n(\widehat{\mathcal M}_{\bm c, 1}(\bm t))|^2)|\Mela{g}{\bm c}(\bm t)|^{-2}d\bm t \leq \sigma n^{-1} \Delta_g(\bm k')=: \Psi^2.
	\end{align*}
	Now for $\tau$ we see that $\Var_{f_Y}(\nu_h(Y_1)) \leq \E_{f_Y}(\nu_h^2(Y_1)) \leq \| f_Y \bm x^{\ushort{2\bm c-\bm 1}}\|_{\infty} \| \nu_{h} \|_{\bm x^{\ushort{\bm 1-2\bm c}}}^2$. Further,
	\begin{align*}
	\| \nu_{h} \|_{\bm x^{\ushort{\bm 1-2\bm c}}}^2 = (2\pi)^{-d} \int_{Q_{\bm k'}}|\Mela{h}{\bm c}(\bm t)|^2|\Mela{g}{\bm c}(\bm t)|^{-2} d\bm t \leq \|\1_{Q_{\bm k'}} \Mela{g}{\bm c}^{-2}\|_{\infty}.
	\end{align*}
	Thus we choose $\tau:= \|f_Y \bm x^{ \ushort{2\bm c-\bm 1}} \|_{\infty} \|\1_{Q_{\bm k'}} \Mela{g}{\bm c}^{-2}\|_{\infty}$. Let us now consider $\psi^2$. We have for any $\bm y \in \pRz^d$,
	\begin{align*}
	|\nu_h(\bm y)|^2 = (2\pi)^{-2d} \left| \int_{Q_{\bm k'}} \bm y^{\ushort{\bm c-\bm 1+i\bm t}} \1_{(0, c_n)}(\bm y^{\ushort{\bm c-\bm 1}}) \frac{\Mela{h}{\bm c}(-\bm t)}{\Mela{g}{\bm c}(\bm t)} d\bm t\right|^2 \leq c_n^2 \Delta_g(\bm k)=:\psi^2,
	\end{align*}
	since $\|h\|_{\bm x^{\ushort{2\bm c-\bm1}}} \leq 1$ and $|\bm y^{\ushort{i\bm t}}|=1$. Applying now the Talagrand inequality we get
	\begin{align*}
	\E_{f_Y}^n((\sup_{h\in \Uz} \overline\nu_h^2 -6\Psi^2)_+) &\leq \frac{C_{f_Y}}{n} \left(\bm k^{\ushort{2\bm \gamma}} \exp(-C_{f_Y,\sigma} \bm k^{\ushort{\bm 1}}) + c_n^2 \exp(-\frac{\sqrt{n\sigma}}{100 c_n})\right)  \leq \frac{C_{f_Y, \sigma}}{n} \left(\bm k^{\ushort{2\bm \gamma}} \exp(-C_{f_Y,\sigma} \bm k^{\ushort{\bm 1}}) + n^{-d} \right)
	\end{align*}
	for the choice $c_n:= \sqrt{n\sigma}/(100\log(n^{d+1}))$. For $\chi \geq 144$ we can conclude that
	\begin{align*}
	\hspace*{-0.7cm}\E_{f_Y}^n(\left( \sup_{h\in \Uz} \nu_{\widehat{\bm k}, 1}(h)^2- \frac{\chi}{24} \sigma\Delta_g(\widehat{\bm k}) n^{-1}\right)_+)  &\leq \sum_{\bm k\in \mathcal K_n}\frac{C_{f_Y, \sigma}}{n} \left(\bm k^{\ushort{2\bm \gamma}} \exp(-C_{f_Y,\sigma} \bm k^{\ushort{\bm 1}}) + n^{-d} \right) \leq C_{f_Y, \sigma, \gamma} n^{-1}
	\end{align*}
	since $|\mathcal K_n| \leq n^{d}$.
	For the second summand in \eqref{eq:two:process} we get for any $\bm k'\in \mathcal K_n$ and $h\in \Uz$,
	\begin{align*}
	|\nu_{\widehat {\bm k}, 2}(h)|^2 &\leq (2\pi)^{-d} \int_{Q_{\widehat{\bm k}}} |\widehat{\mathcal M}_{\bm c, 2}(\bm t)- \E_{f_Y}^n(\widehat{\mathcal M}_{\bm c, 2}(\bm t))|^2|\Mela{g}{\bm c}(\bm t)|^{-2} d\bm t \\
	& \leq\sum_{\bm k'\in \mathcal K_n} (2\pi)^{-d} \int_{Q_{\bm{k}'}} |\widehat{\mathcal M}_{\bm c, 2}(\bm t)- \E_{f_Y}^n(\widehat{\mathcal M}_{\bm c, 2}(\bm t))|^2|\Mela{g}{\bm c}(\bm t)|^{-2} d\bm t.
	\end{align*}
	Thus we have for any $u>0$
	\begin{align*}
	\hspace*{-1cm}\E_{f_Y}^n(\sup_{h\in \Uz}\nu_{\widehat{\bm k},2}(h)^2)) \leq \sum_{\bm k' \in \mathcal K_n}\Delta_g(\bm{k'}) n^{-1} \E_{f_Y}^1(\bm Y_1^{\ushort{2\bm c-\bm 2}} \1_{[c_n, \infty)}(\bm Y_1^{\ushort{\bm c-\bm 1}}))
	\leq \frac{C_g|\mathcal K_n|}{c_n^u} \E_{f_Y}(\bm Y_1^{\ushort{(2+u)(\bm c-\bm 1)}}).
	\end{align*}
	Now under assumption \eqref{eq:ass:g1} we have $|\mathcal K_n| \leq |\{\bm k \in \mathbb N^d: \bm k^{\ushort{2\bm \gamma+ \bm 1}} \leq c_g n\}| \leq C_g n \log(n)^{d-1}$, compare \cite{Dussap2022}.
	Now choosing $u=5$ implies 
	$$\E_{f_Y}^n( \sup_{h\in \Uz}\nu_{\widehat{\bm k},2}(h)^2)) \leq C_{g,\sigma} \E_{f_Y}(\bm Y_1^{\ushort{7(\bm c-\bm 1)}}) n^{-1}.$$
	To finish the proof we still need to show (ii). To do so, we define the event $\Omega:=\{|\widehat \sigma-\sigma| \leq \sigma/2\}$. On $\Omega$ holds $\sigma \leq 2\widehat \sigma \leq 3 \sigma $ and we deduce
	$$ (\mathrm{pen}(\widehat{\bm k})-2\widehat{\mathrm{pen}}(\widehat{\bm k}))_+ \leq \chi (\sigma-2\widehat\sigma)_+ \mathds 1_{\Omega^c}$$
	since $\widehat{\bm k} \in \mathcal K_n$. Therefore we get by application of the Cauchy-Schwartz inequality and the Markow inequality 
	$$ \E((\mathrm{pen}(\widehat{\bm k})-2\widehat{\mathrm{pen}}(\widehat{\bm k}))_+) \leq C(\chi, \sigma) \Var(\widehat\sigma) = C(\chi, \sigma, \E(\bm Y_1^{\ushort{4(\bm c- \bm 1)}})) n^{-1}.$$
\end{proof}

\section*{References}

%\cite{andersson,anderson1993totally,lauritzen1996graphical,gratzer2002general} or \cite{drton2008iterative}. At the end of this template "trial.bib" is included and below there is also an example of how to include references in an alternative way.
%\section*{References}

% To ensure accuracy, get them from MathSciNet whenever possible. Typeset them with BibTeX using JMVA's style file, \texttt{myjmva.bst}.
%\bibliography{}
\bibliographystyle{myjmva}
%\begin{thebibliography}
\bibliography{AISCEMME}

\begin{thebibliography}{18}
\expandafter\ifx\csname natexlab\endcsname\relax\def\natexlab#1{#1}\fi
\providecommand{\bibinfo}[2]{#2}
\ifx\xfnm\relax \def\xfnm[#1]{\unskip,\space#1}\fi
%Type = Article
\bibitem[{Belomestny and Goldenshluger(2020)}]{BelomestnyGoldenshluger2020}
\bibinfo{author}{D.~Belomestny}, \bibinfo{author}{A.~Goldenshluger},
  \bibinfo{title}{Nonparametric density estimation from observations with
  multiplicative measurement errors}, \bibinfo{journal}{Ann. Inst. Henri
  Poincar\'{e} Probab. Stat.} \bibinfo{volume}{56} (\bibinfo{year}{2020})
  \bibinfo{pages}{36--67}.
%Type = Article
\bibitem[{Brenner~Miguel(2021)}]{Brenner-Miguel2021}
\bibinfo{author}{S.~Brenner~Miguel}, \bibinfo{title}{Multiplicative
  deconvolution estimator based on a ridge approach}, \bibinfo{journal}{arXiv
  preprint arXiv:2108.01523}  (\bibinfo{year}{2021}).
%Type = Article
\bibitem[{Brenner~Miguel et~al.(2021)Brenner~Miguel, Comte and
  Johannes}]{Brenner-MiguelComteJohannes2020}
\bibinfo{author}{S.~Brenner~Miguel}, \bibinfo{author}{F.~Comte},
  \bibinfo{author}{J.~Johannes}, \bibinfo{title}{{Spectral cut-off
  regularisation for density estimation under multiplicative measurement
  errors}}, \bibinfo{journal}{Electronic Journal of Statistics}
  \bibinfo{volume}{15} (\bibinfo{year}{2021}) \bibinfo{pages}{3551 -- 3573}.
%Type = Article
\bibitem[{Comte and Dion(2016)}]{ComteDion2016}
\bibinfo{author}{F.~Comte}, \bibinfo{author}{C.~Dion},
  \bibinfo{title}{Nonparametric estimation in a multiplicative censoring model
  with symmetric noise}, \bibinfo{journal}{J. Nonparametr. Stat.}
  \bibinfo{volume}{28} (\bibinfo{year}{2016}) \bibinfo{pages}{768--801}.
%Type = Article
\bibitem[{Comte and Lacour(2013)}]{ComteLacour2013}
\bibinfo{author}{F.~Comte}, \bibinfo{author}{C.~Lacour},
  \bibinfo{title}{Anisotropic adaptive kernel deconvolution},
  \bibinfo{journal}{Ann. Inst. Henri Poincar\'{e} Probab. Stat.}
  \bibinfo{volume}{49} (\bibinfo{year}{2013}) \bibinfo{pages}{569--609}.
%Type = Article
\bibitem[{Dussap(2021)}]{Dussap2021}
\bibinfo{author}{F.~Dussap}, \bibinfo{title}{Anisotropic multivariate
  deconvolution using projection on the {L}aguerre basis}, \bibinfo{journal}{J.
  Statist. Plann. Inference} \bibinfo{volume}{215} (\bibinfo{year}{2021})
  \bibinfo{pages}{23--46}.
%Type = Article
\bibitem[{Dussap(2022)}]{Dussap2022}
\bibinfo{author}{F.~Dussap}, \bibinfo{title}{Nonparametric multiple regression
  by projection on non-compactly supported bases}  (\bibinfo{year}{2022}).
%Type = Book
\bibitem[{Engl et~al.(1996)Engl, Hanke and Neubauer}]{EnglHankeNeubauer1996}
\bibinfo{author}{H.~W. Engl}, \bibinfo{author}{M.~Hanke},
  \bibinfo{author}{A.~Neubauer}, \bibinfo{title}{Regularization of inverse
  problems}, volume \bibinfo{volume}{375} of \text{\bibinfo{series}{Mathematics
  and its Applications}}, \bibinfo{publisher}{Kluwer Academic Publishers Group,
  Dordrecht}, \bibinfo{year}{1996}.
%Type = Article
\bibitem[{Fan(1991)}]{Fan1991}
\bibinfo{author}{J.~Fan}, \bibinfo{title}{On the optimal rates of convergence
  for nonparametric deconvolution problems}, \bibinfo{journal}{Ann. Statist.}
  \bibinfo{volume}{19} (\bibinfo{year}{1991}) \bibinfo{pages}{1257--1272}.
%Type = Article
\bibitem[{Klein and Rio(2005)}]{KleinRio2005}
\bibinfo{author}{T.~Klein}, \bibinfo{author}{E.~Rio},
  \bibinfo{title}{Concentration around the mean for maxima of empirical
  processes}, \bibinfo{journal}{Ann. Probab.} \bibinfo{volume}{33}
  (\bibinfo{year}{2005}) \bibinfo{pages}{1060--1077}.
%Type = Article
\bibitem[{Lepski and Willer(2019)}]{LepskiWiller2019}
\bibinfo{author}{O.~V. Lepski}, \bibinfo{author}{T.~Willer},
  \bibinfo{title}{Oracle inequalities and adaptive estimation in the
  convolution structure density model}, \bibinfo{journal}{Ann. Statist.}
  \bibinfo{volume}{47} (\bibinfo{year}{2019}) \bibinfo{pages}{233--287}.
%Type = Article
\bibitem[{Liu et~al.(2013)Liu, Xiao and Wu}]{LiuXiaoWu2013}
\bibinfo{author}{W.~Liu}, \bibinfo{author}{H.~Xiao}, \bibinfo{author}{W.~B.
  Wu}, \bibinfo{title}{Probability and moment inequalities under dependence},
  \bibinfo{journal}{Statist. Sinica} \bibinfo{volume}{23}
  (\bibinfo{year}{2013}) \bibinfo{pages}{1257--1272}.
%Type = Book
\bibitem[{Meister(2009)}]{Meister2009}
\bibinfo{author}{A.~Meister}, \bibinfo{title}{Deconvolution problems in
  nonparametric statistics}, volume \bibinfo{volume}{193} of
  \text{\bibinfo{series}{Lecture Notes in Statistics}},
  \bibinfo{publisher}{Springer-Verlag, Berlin}, \bibinfo{year}{2009}.
%Type = Article
\bibitem[{Nagaev(1979)}]{Nagaev1979}
\bibinfo{author}{S.~V. Nagaev}, \bibinfo{title}{Large deviations of sums of
  independent random variables}, \bibinfo{journal}{Ann. Probab.}
  \bibinfo{volume}{7} (\bibinfo{year}{1979}) \bibinfo{pages}{745--789}.
%Type = Book
\bibitem[{Paris and Kaminski(2001)}]{ParisKaminski2001}
\bibinfo{author}{R.~B. Paris}, \bibinfo{author}{D.~Kaminski},
  \bibinfo{title}{Asymptotics and {M}ellin-{B}arnes integrals},
  volume~\bibinfo{volume}{85} of \text{\bibinfo{series}{Encyclopedia of
  Mathematics and its Applications}}, \bibinfo{publisher}{Cambridge University
  Press, Cambridge}, \bibinfo{year}{2001}.
%Type = Article
\bibitem[{Rebelles(2016)}]{Rebelles2016}
\bibinfo{author}{G.~Rebelles}, \bibinfo{title}{Structural adaptive
  deconvolution under {$\Bbb{L}_p$}-losses}, \bibinfo{journal}{Math. Methods
  Statist.} \bibinfo{volume}{25} (\bibinfo{year}{2016})
  \bibinfo{pages}{26--53}.
%Type = Article
\bibitem[{Talagrand(1996)}]{Talagrand1996}
\bibinfo{author}{M.~Talagrand}, \bibinfo{title}{New concentration inequalities
  in product spaces}, \bibinfo{journal}{Invent. Math.} \bibinfo{volume}{126}
  (\bibinfo{year}{1996}) \bibinfo{pages}{505--563}.
%Type = Book
\bibitem[{Tsybakov(2009)}]{Tsybakov2009}
\bibinfo{author}{A.~B. Tsybakov}, \bibinfo{title}{Introduction to nonparametric
  estimation}, Springer Series in Statistics, \bibinfo{publisher}{Springer, New
  York}, \bibinfo{year}{2009}. \bibinfo{note}{Revised and extended from the
  2004 French original, Translated by Vladimir Zaiats}.

\end{thebibliography}
%\end{document}

%\bibliographystyle{myjmva}
%\section*{}
%or one can use (see the template for details)
%\begin{thebibliography}{99}

%\bibitem[{Agresti(2013)}]{Agresti13}
%\bibinfo{author}{A.~Agresti}, \bibinfo{title}{{Categorical Data Analysis}},
%  \bibinfo{publisher}{Wiley, Hoboken}, \bibinfo{year}{2013}.
%Type = Article
%\bibitem[{Aitchison and Silvey(1958)}]{aitchison1958maximum}
%\bibinfo{author}{J.~Aitchison}, \bibinfo{author}{S.D.~Silvey},
%  \bibinfo{title}{Maximum-likelihood estimation of parameters subject to
%  restraints}, \bibinfo{journal}{The Annals of Mathematical Statistics}
%  \bibinfo{volume}{29} (\bibinfo{year}{1958}) \bibinfo{pages}{813--828}.
% \bibitem{Balsu} A. Balsubramani, S. Dasgupta, Y. Freund, \newblock The fast
%convergence of incremental PCA, \newblock Advances in Neural Information
%Processing Systems 26 (2013) 3174--3182.

%\end{thebibliography}

\end{document}